\documentclass[11pt]{article}
\usepackage{epsfig}

\def\be{\begin{equation}}
\def\ee{\end{equation}}
\def\bea{\begin{eqnarray}}
\def\eea{\end{eqnarray}}
\def\bes{\begin{eqnarray*}}
\def\ees{\end{eqnarray*}}

\def\nn{\nonumber}
\def\<{\langle}
\def\>{\rangle}
\def\lb{\label}
\def\bs{\setminus}

\def\R{{\bf R}}
\def\C{{\bf C}}
\def\Z{{\bf Z}}
\def\N{{\bf N}}

\def\Q{{\bf Q}}
\def\T{{\bf T}}

\def\aa{{\alpha}}
\def\bb{{\beta}}
\def\ga{{\gamma}}

\def\th{{\theta}}

\def\Om{{\Omega}}
\def\ep{{\epsilon}}
\def\lm{{\lambda}}
\def\Lm{{\Lambda}}

\def\Sg{{\Sigma}}
\def\vf{{\varphi}}

\def\H{{\cal H}}

\def\T{{\cal T}}

\def\Nn{{\cal N}}

\def\mul{{\rm mul}}

\def\Sp{{\rm Sp}}
\def\mod{{\rm mod}}
\def\crit{{\rm crit}}

\def\dm{{\rm \diamond}}

\def\wtd#1{\widetilde{#1}}
\def\hb{\vrule height0.18cm width0.14cm $\,$}
\def\mapright#1{\smash{\mathop{\longrightarrow}\limits^{#1}}}

\title{Multiplicity and ellipticity of closed characteristics on compact star-shaped hypersurfaces in $\R^{2n}$}
\author{
Huagui Duan$^{1}$,\thanks{Partially supported by NSFC (Nos. 11131004, 11471169),
LPMC of MOE of China and Nankai University. E-mail: duanhg@nankai.edu.cn.}
\quad Hui Liu$^{2}$\thanks{Partially supported by NSFC (Nos.11401555, 11371339), Anhui Provincial Natural Science Foundation (No. 1608085QA01).
E-mail: huiliu@ustc.edu.cn.}  \\ \\
$^{1}$ School of Mathematical Sciences and LPMC, Nankai University,\\Tianjin 300071, People's Republic of China\\
$^{2}$ Key Laboratory of Wu Wen-Tsun Mathematics, Chinese Academy of Sciences,\\
School of Mathematical Sciences, University of Science and Technology of China,
\\Hefei, Anhui 230026, People's Republic of China\\}

\begin{document}
\date{December 21, 2015}

\maketitle
\begin{abstract}
{\it In this paper, we firstly generalize some theories developed by I. Ekeland and H. Hofer in \cite{EkH1} for
closed characteristics on compact convex hypersurfaces in $\R^{2n}$ to star-shaped hypersurfaces. As applications
we use Ekeland-Hofer theory and index iteration theory to prove that if a compact star-shaped hypersuface in ${\bf R}^4$
satisfying some suitable pinching condition carries exactly two geometrically distinct closed characteristics, then both of them must be elliptic.
We also conclude that the theory developed by Y. Long and C. Zhu in \cite{LoZ1} still holds for dynamically convex star-shaped hypersurfaces,
and combining it with the results in \cite{WHL1}, \cite{LLW1}, \cite{Wan3}, we obtain that there exist
at least $n$ closed characteristics on every dynamically convex star-shaped hypersurface in $\R^{2n}$ for $n=3, 4$.  }
\end{abstract}


\renewcommand{\theequation}{\thesection.\arabic{equation}}
\renewcommand{\thefigure}{\thesection.\arabic{figure}}

\setcounter{equation}{0}
\section{Introduction and main result}

Let $\Sigma$ be a $C^3$ compact hypersurface in $\R^{2n}$ strictly star-shaped with respect to the origin, i.e.,
the tangent hyperplane at any $x\in\Sigma$ does not intersect the origin. We denote the set of all such
hypersurfaces by $\H_{st}(2n)$, and denote by $\H_{con}(2n)$ the subset of $\H_{st}(2n)$ which consists of all
strictly convex hypersurfaces. We consider closed characteristics $(\tau, y)$ on $\Sigma$, which are solutions
of the following problem
\be
\left\{\matrix{\dot{y}=JN_\Sigma(y), \cr
               y(\tau)=y(0), \cr }\right. \lb{1.1}\ee
where $J=\left(\matrix{0 &-I_n\cr
        I_n  & 0\cr}\right)$, $I_n$ is the identity matrix in $\R^n$, $\tau>0$, $N_\Sigma(y)$ is the outward
normal vector of $\Sigma$ at $y$ normalized by the condition $N_\Sigma(y)\cdot y=1$. Here $a\cdot b$ denotes
the standard inner product of $a, b\in\R^{2n}$. A closed characteristic $(\tau, y)$ is {\it prime}, if $\tau$
is the minimal period of $y$. Two closed characteristics $(\tau, y)$ and $(\sigma, z)$ are {\it geometrically
distinct}, if $y(\R)\not= z(\R)$. We denote by $\T(\Sigma)$ the set of geometrically distinct
closed characteristics $(\tau, y)$ on $\Sigma\in\mathcal{H}_{st}(2n)$. A closed characteristic
$(\tau,y)$ is {\it non-degenerate} if $1$ is a Floquet multiplier of $y$ of precisely algebraic multiplicity
$2$; {\it hyperbolic} if $1$ is a double Floquet multiplier of it and all the other Floquet multipliers
are not on ${\bf U}=\{z\in {\bf C}\mid |z|=1\}$, i.e., the unit circle in the complex plane; {\it elliptic}
if all the Floquet multipliers of $y$ are on ${\bf U}$. We call a $\Sigma\in \mathcal{H}(2n)$ {\it non-degenerate} if
all the closed characteristics on $\Sigma$ together with all of their iterations are non-degenerate.

Fix a constant $\alpha$ satisfying $1<\alpha<2$ and define the Hamiltonian function
$H_\alpha:\R^{2n}\to [0,+\infty)$ by \be  H_\alpha(x) = j(x)^{\alpha}, \qquad \forall x\in\R^{2n},  \lb{1.2}\ee
where $j$ is the gauge function of $\Sigma$, i.e., $j(x)=\lm$ if $x=\lm y$ for
some $\lm>0$ and $y\in\Sigma$ when $x\in\R^{2n}\bs\{0\}$, and
$j(0)=0$. Then $H_{\alpha}\in C^1(\R^{2n},\R)\cap C^3(\R^{2n}\bs\{0\},\R)$ and
$\Sigma=H_{\alpha}^{-1}(1)$. It is well-known that the problem
(\ref{1.1}) is equivalent to the following given energy problem of
the Hamiltonian system \bea \left\{\begin{array}{ll} \dot{y}(t) =
JH_{\alpha}^{\prime}(y(t)), \quad H_{\alpha}(y(t)) =1, \qquad
\forall t\in\R, \cr y(\tau) = y(0). \end{array}\right. \lb{1.3}\eea
 Denote by
$\T(\Sigma,\alpha)$ the set of all geometrically distinct solutions
$(\tau,y)$ of the problem (\ref{1.3}). Note that elements in
$\T(\Sigma)$ and $\T(\Sigma,\alpha)$ are one to
one correspondent to each other.

The study on closed characteristics in the global sense started in 1978, when the existence of at least one
closed characteristic was first established on any $\Sg\in\H_{st}(2n)$ by P. Rabinowitz in \cite{Rab1}
and on any $\Sg\in\H_{con}(2n)$ by A. Weinstein in \cite{Wei1} independently, since then the existence of
multiple closed characteristics on $\Sg\in\H_{con}(2n)$ has been deeply studied by many mathematicians, for
example, studies in \cite{EkL1}, \cite{EkH1}, \cite{Szu1}, \cite{HWZ1}, \cite{LoZ1}, \cite{WHL1}, \cite{Wan2} and \cite{Wan3} for
convex hypersurfaces. For the star-shaped hypersurfaces, in \cite{Gir1} of 1984 and \cite{BLMR} of 1985,
$\;^{\#}\T(\Sg)\ge n$ for $\Sg\in\H_{st}(2n)$ was proved under some pinching conditions. In \cite{Vit1} of 1989,
C. Viterbo proved a generic existence result for infinitely many closed characteristics on star-shaped
hypersurfaces. In \cite{HuL1} of 2002, X. Hu and Y. Long proved that $\;^{\#}\T(\Sg)\ge 2$ for
$\Sg\in \H_{st}(2n)$ on which all the closed characteristics and their iterates are non-degenerate.
In \cite{HWZ2} of 2003, H. Hofer, K. Wysocki, and E. Zehnder proved any non-degenerate compact
star-shaped hypersurface has either two or infinitely closed characteristics,
provided that all stable and unstable manifolds of the hyperbolic closed characteristics intersect transversally.
Recently $\;^{\#}\T(\Sg)\ge 2$ was first proved for every $\Sg\in \H_{st}(4)$
by D. Cristofaro-Gardiner and M. Hutchings in \cite{CGH1} without any pinching or non-degeneracy conditions.
Different proofs of this result can also be found in \cite{GHHM}, \cite{LLo1} and \cite{GiG1}.

I. Ekeland and H. Hofer in \cite{EkH1} provided a close relationship
between the set of Maslov-type indices of closed characteristics and the set of even
positive integers, which is the core in studying the multiplicity and ellipticity of the
closed characteristics on compact convex hypersurfaces (cf. \cite{LoZ1}). Our main goal in this paper is to generalize
the theory of Ekeland-Hofer to compact star-shaped hypersurfaces and as its applications we give some multiplicity and stability
results of closed characteristics on compact star-shaped hypersurfaces.

For the stability of closed characteristics on $\Sigma \in\H_{st}(2n)$
we refer the readers to \cite{LiL1} and \cite{LLo2}. Specially, in \cite{LLo2},
H. Liu and Y. Long proved that $\Sg\in\H_{st}(4)$ and $\,^{\#}\T(\Sg)=2$ imply that both of the closed
characteristics must be elliptic provided that $\Sg$ is symmetric with respect to the origin.

Let $n(y)$ be the unit outward normal vector of $\Sigma$ at $y$ and
$d(y):=n(y)\cdot y$, i.e., the distance between the origin of $\R^{2n}$ and the tangent hyperplane to $\Sigma$
at $y$, then $d(y)>0$ for all $y\in \Sigma$ since $\Sigma$ is strictly star-shaped. Let $d=\min{\{d(y): y\in \Sigma\}}$ and
$R=\max{\{|y|: y\in \Sigma\}}$. In this paper, we prove, under suitable pinching condition, the symmetric condition in
Theorem 1.4 of \cite{LLo2} can be dropped, i.e., the following theorem holds.

\medskip

{\bf Theorem 1.1.} {\it  Suppose that $\Sg\in \mathcal{H}_{st}(4)$ satisfy $\,^{\#}\T(\Sg)=2$ and $R^2< 2d^2$. Then both of
the closed characteristics are elliptic.}

\medskip

{\bf Remark 1.2.} Note that the pinching condition on $\Sg$ in Theorem 1.1 is only used to get a
contradiction in the study of the Subcase 1.2 of Case 1 in the proof of Theorem 1.1.

\medskip

In Definition 1.2 of \cite{HWZ1}, An interesting class of contact forms on $S^3$ which are called dynamically
convex contact forms was introduced. Similarly, we give the following definition:

\medskip

{\bf Definition 1.3.} {\it $\Sg\in\H_{st}(2n)$ is called dynamically
convex if any closed characteristic $(\tau, y)$ on $\Sigma$ has its Maslov-type index not less than $n$.}

\medskip

Note that from the proof of Theorem 3.4 of \cite{HWZ1}, for $n=2$, the above definition coincides
with that of \cite{HWZ1}. Also from the Remark before Definition 3.6 of \cite{HWZ1} and
Corollary 1.2 of \cite{LoZ1}, we know any $\Sg\in\H_{con}(2n)$ is dynamically convex. As mentioned before
Definition 1.2 of \cite{HWZ1}, ``strictly convex" is not a symplectically invariant
concept, thus a dynamically convex $\Sg\in\H_{st}(2n)$ need not to be convex.

In this paper, we also prove the main results of \cite{LoZ1}, \cite{WHL1}, \cite{Wan2}, \cite{Wan3} and \cite{HuO} hold for
dynamically convex star-shaped hypersurfaces which cover the works of these literatures.
Specially, we have:

\medskip

{\bf Theorem 1.4.} {\it Let $\Sg\in\H_{st}(2n)$ be dynamically convex. Then
$^\#\T(\Sg)\geq [\frac{n+1}{2}]+1$. If $\Sigma$ is nondegenerate, then $^\#\T(\Sg)\geq n$.
If $^\#\T(\Sg)<+\infty$, then there exists at least two elliptic closed characteristic on $\Sigma$,
and at least $[n/2]$ closed characteristics possessing irrational mean indices.  }

\medskip

{\bf Theorem 1.5.} {\it Let $\Sg\in\H_{st}(2n)$ be dynamically convex. Then
$^\#\T(\Sg)\geq n$ for $n=3, 4$.}

\medskip

{\bf Remark 1.6.} Note that J. Gutt and J. Kang in \cite{GuK} proved that if $\Sigma\in\mathcal{H}_{st}(2n)$ is non-degenerate and dynamically convex,
then there exist at least $n$ closed characteristics on such $\Sg$, whose iterates' Conley-Zehnder indices possess the same parity. Note
that their index definition is slightly different from ours. Also recently, M. Abreu and L. Macarini
in \cite{AbM} gave a sharp lower bound for the number of geometrically distinct contractible periodic
orbits of non-degenerate dynamically convex Reeb flows on prequantizations of symplectic manifolds that are not
aspherical, which implies results of \cite{GuK} (cf. Corollary 2.9 of \cite{AbM}).
We also mention that very recently, Y. Long, W. Wang and the authors in \cite{DLLW} proved some sharp multiplicity results
for non-degenerate star-shaped hypersurfaces under some index conditions, which are weaker than the convex or dynamically convex case.
Our Theorem 1.4 and Theorem 1.5 give new multiplicity and stability results for the degenerate, dynamically convex, star-shaped hypersurfaces.

\medskip

This paper is arranged as follows. In Sections 2 and 3, following the frame works of \cite{Vit1}, \cite{LLW1} and \cite{Eke1},
we establish a variational structure for closed characteristics on
star-shaped hypersurfaces and prove some theories in \cite{EkH1} hold for star-shaped
case, we omit most of the details of the proofs of the theories below
and only point out differences from \cite{Vit1}, \cite{LLW1} and \cite{Eke1} when necessary.
In Section 3.2, we further study the critical values obtained in Section 3.1 when the star-shaped hypersurface is suitably pinched,
which we will use to prove Theorem 1.1 in Section 4.  In Section 4, we also explain how to get Theorems 1.4 and 1.5 as another application of Ekeland-Hofer theory.
In Section 5 (an appendix), we briefly review the equivariant Morse theory and the resonance identities for closed characteristics on
compact star-shaped hypersurfaces in $\R^{2n}$ developed in \cite{LLW1}, which are used in the proof of Theorem 1.1.

In this paper, let $\N$, $\N_0$, $\Z$, $\Q$, $\R$, $\C$ and $\R^+$ denote the sets of natural integers,
non-negative integers, integers, rational numbers, real numbers, complex numbers and positive real numbers respectively.
We define the function $[a]=\max{\{k\in {\bf Z}\mid k\leq a\}}$, $\{a\}=a-[a]$ , and
$E(a)=\min{\{k\in{\bf Z}\mid k\geq a\}}$.
Denote by $a\cdot b$ and $|a|$ the standard
inner product and norm in $\R^{2n}$. Denote by $\langle\cdot,\cdot\rangle$ and $\|\cdot\|$
the standard $L^2$ inner product and $L^2$ norm. For an $S^1$-space $X$, we denote by
$X_{S^1}$ the homotopy quotient of $X$ by $S^1$, i.e., $X_{S^1}=S^\infty\times_{S^1}X$,
where $S^\infty$ is the unit sphere in an infinite dimensional complex Hilbert space.
In this paper we use $\Q$ coefficients for all homological and cohomological modules. By $t\to a^+$, we
mean $t>a$ and $t\to a$.

\setcounter{equation}{0}
\section{A variational structure for closed characteristics on compact star-shaped hypersurfaces}

In Section 2 and 3, we fix a $\Sg\in\H_{st}(2n)$.

As in Sections V.2 and V.3 of \cite{Eke1}, we consider the following fixed period problem

\be \left\{\matrix{\dot{x}(t) &=& JH_\alpha^\prime(x(t)), \cr
                   x(0) &=& x(1).\qquad    \cr }\right.   \lb{2.1}\ee
Then solutions of (\ref{2.1}) are $x\equiv 0$ and $x=\tau^{-\frac{1}{2-\alpha}} y(\tau t)$, where $(\tau, y)$ is a solution of (\ref{1.3}).

For technical reasons (to get Proposition 2.5 below), we need to further modify the Hamiltonian, more precisely,
we follow Page 624 of \cite{Vit2}, and let $\ep$ satisfy $\ep<2\pi$, we can construct a function $H$, which coincides
with $H_\alpha$ on $U_A=\{x\mid H_\alpha(x)\leq A\}$ for some large $A$, and with $\frac{1}{2}\ep|x|^2$ outside
some large ball, such that $\nabla H(x)$ does not vanish and $H^{\prime\prime}(x)<\ep$ outside $U_A$. As in
Proposition 2.7 of \cite{Vit2}, we have the following result.

{\bf Proposition 2.1.} {\it For small $\ep$, there exists a function $H$ on ${\bf R}^{2n}$ such that $H$ is $C^1$ on ${\bf R}^{2n}$,
and $C^3$ on ${\bf R}^{2n}\bs\{0\}$, $H=H_\alpha$ in $U_A$, and $H(x)=\frac{1}{2}\ep|x|^2$ for $|x|$ large,
and the solutions of the fixed period system
\be   \left\{\matrix{\dot{x}(t) &=& JH^\prime(x(t)), \cr
                   x(0) &=& x(1),\qquad         \cr }\right.  \lb{2.2}\ee
are the same with those of (\ref{2.1}), i.e., the solutions of (\ref{2.2}) are $x\equiv 0$
and $x=\tau^{-\frac{1}{2-\alpha}} y(\tau t)$, where $(\tau, y)$ is a solution of (\ref{1.3}).}

Note that the condition (2.2) of Lemma 2.2 of \cite{Vit2} is only used to get Theorem 7.1 of \cite{Vit2}, so the
other statements in \cite{Vit2} also hold for our choice of the Hamiltonian function.

As in \cite{BLMR} (cf. Section 3 of \cite{Vit2}), we can choose some
large constant $K$ such that
\be H_{K}(x) = H(x)+\frac{1}{2}K|x|^2   \lb{2.3}\ee
is a strictly convex function, that is,
\be (\nabla H_{K}(x)-\nabla H_{K}(y), x-y) \geq \frac{\ep}{2}|x-y|^2,  \lb{2.4}\ee
for all $x, y\in {\bf R}^{2n}$, and some positive $\ep$. Let $H_{K}^*$ be the Fenchel dual of $H_{K}$
defined by
\be  H_{K}^\ast (y) = \sup\{x\cdot y-H_{K}(x)\;|\; x\in \R^{2n}\}.   \lb{2.5}\ee
The dual action functional on $X=W^{1, 2}({\bf R}/{{\bf Z}}, {\bf R}^{2n})$ is defined by
\be F_{K}(x) = \int_0^1{\left[\frac{1}{2}(J\dot{x}-K x,x)+H_{K}^*(-J\dot{x}+K x)\right]dt}.  \lb{2.6}\ee
Then we have

{\bf Lemma 2.2.} (cf. Proposition 3.4 of \cite{Vit2}) {\it Assume $K\not\in 2\pi{\bf Z}$, then $x$ is a
critical point of $F_{K}$ if and only if it is a solution of (\ref{2.2}).}

As is well known, when $K\notin 2\pi{\bf Z}$, the map $x\mapsto -J\dot{x}+Kx$ is a Hilbert space isomorphism between
$X=W^{1, 2}({\bf R}/{{\bf Z}}; {\bf R}^{2n})$ and $E=L^{2}({\bf R}/{\bf Z},{\bf R}^{2n})$. We denote its inverse
by $M_K$ and the functional
\be \Psi_{K}(u)=\int_0^1{\left[-\frac{1}{2}(M_{K}u, u)+H_{K}^*(u)\right]dt}, \qquad \forall\,u\in E. \lb{2.7}\ee
Then $x\in X$ is a critical point of $F_{K}$ if and only if $u=-J\dot{x}+Kx$ is a critical point of $\Psi_{K}$.
We have a natural $S^1$-action on $X$ or $E$ defined by
\be  \theta\cdot u(t)=u(\theta+t),\quad\forall\, \theta\in S^1, \, t\in\R.  \lb{2.8}\ee
Clearly both of $F_{K}$ and $\Psi_{K}$ are $S^1$-invariant. For any $\kappa\in\R$, we denote by
\bea
\Psi_{K}^\kappa = \{u\in L^{2}({\bf R}/{T {\bf Z}}; {\bf R}^{2n})\;|\;\Psi_{K}(u)\le\kappa\}.\nn\eea
Obviously, this level set is also $S^1$-invariant.

{\bf Definition 2.3.} (cf. p.628 of \cite{Vit2}) {\it Suppose $u$ is a nonzero critical point of $\Psi_{K}$.
Then the formal Hessian of $\Psi_{K}$ at $u$ is defined by
\be Q_{K}(v)=\int_0^1(-M_K v\cdot v+H_{K}^{*\prime\prime}(u)v\cdot v)dt,  \lb{2.9}\ee
which defines an orthogonal splitting $E=E_-\oplus E_0\oplus E_+$ into negative, zero and positive subspaces.
The index and nullity of $u$ are defined by $i_K(u)=\dim E_-$ and $\nu_K(u)=\dim E_0$ respectively. }

Similarly, we define the index and nullity of $x=M_Ku$ for $F_{K}$, which are denoted by $i_K(x)$ and
$\nu_K(x)$ respectively. Then we have
\be  i_K(u)=i_K(x),\quad \nu_K(u)=\nu_K(x),  \lb{2.10}\ee
which follow from the definitions (\ref{2.6}) and (\ref{2.7}). The following important formula was proved in
Lemma 6.4 of \cite{Vit2}:
\be  i_K(x) = 2n([K/{2\pi}]+1)+i^v(x) \equiv d(K)+i^v(x),   \lb{2.11}\ee
where the Viterbo index $i^v(x)$ does not depend on K, but only on $H$.

By the proof of Proposition 2 of \cite{Vit1}, we have that $v\in E$ belongs to the null space of $Q_{K}$
if and only if $z=M_K v$ is a solution of the linearized system
\be  \dot{z}(t) = JH^{\prime\prime}(x(t))z(t).  \lb{2.12}\ee
Thus the nullity in (\ref{2.10}) is independent of $K$, which we denote by $\nu^v(x)\equiv \nu_K(u)= \nu_K(x)$.

Since $x$ is a solution of (\ref{2.1}) corresponding to a solution $(\tau, y)$ of (\ref{1.3}), we also denote
$i^v(x)$ and $\nu^v(x)$ by $i^v(y)$ and $\nu^v(y)$ respectively, and define $i(y):=i^v(y)$ and $\nu(y):=\nu^v(y)$. By Theorem 2.1 of \cite{HuL1}, we have:

{\bf Lemma 2.4.} {\it Suppose $\Sg\in \H_{st}(2n)$ and
$(\tau,y)\in \T(\Sigma)$. Then we have
\bea i(y^m)=i(y, m)-n,\quad \nu(y^m)=\nu(y, m),
       \qquad \forall m\in\N, \nn\eea
where $i(y, m)$ and $\nu(y, m)$ are the Maslov-type index and nullity of $(m\tau,y)$.}

By Propositions 3.9, 4.1 of \cite{Vit2} and the same proof of Proposition 2.12 of \cite{LLW1}, we have:

{\bf Proposition 2.5.} {\it $\Psi_{K}$ satisfies the Palais-Smale condition on $E$, and $F_{K}$ satisfies
the Palais-Smale condition on $X$, when $K\notin 2\pi{\bf Z}$.}

\setcounter{equation}{0}
\section{Fadell-Rabinowitz index theory for closed characteristics on star-shaped hypersurfaces}
\subsection{Critical values in the free case}
Recall that for a principal $U(1)$-bundle $E\to B$, the Fadell-Rabinowitz index
(cf. \cite{FaR1}) of $E$ is defined to be $\sup\{k\;|\, c_1(E)^{k-1}\not= 0\}$,
where $c_1(E)\in H^2(B,\Q)$ is the first rational Chern class. For a $U(1)$-space,
i.e., a topological space $X$ with a $U(1)$-action, the Fadell-Rabinowitz index is
defined to be the index of the bundle $X\times S^{\infty}\to X\times_{U(1)}S^{\infty}$,
where $S^{\infty}\to CP^{\infty}$ is the universal $U(1)$-bundle.

For any $\kappa\in\R$, we denote by
\be \Psi_K^{\kappa-}=\{u\in L^{2}({\bf R}/{\bf Z},{\bf R}^{2n})\;|\;
             \Psi_K(u)<\kappa\}. \lb{3.1}\ee
Then as in P.218 of \cite{Eke1}, we define
\be c_i=\inf\{\delta\in\R\;|\: \hat I(\Psi_K^{\delta-})\ge i\}, i\in\N.\lb{3.2}\ee
where $\hat I$ is the Fadell-Rabinowitz index given above.

For $i\geq d(K)/2+1$,
where $d(K)=2n([K/{2\pi}]+1)$, $c_i$ is well defined.
In fact, by Proposition 5.7 of \cite{Vit2}, there exists constant $c$ such that $\Psi_K^{c}$ is $S^1$-equivariant
homotopy equivalent with a $(d(K)-1)$ dimensional sphere. Then $\hat I(\Psi_K^{c})=d(K)/2$.
Hence for $i\geq d(K)/2+1$, $c_i\geq c$ is well defined.

Then similar to Proposition 3 in P.218 of \cite{Eke1}, we have

{\bf Proposition 3.1.} {\it For $i\geq d(K)/2+1$, $c_i$ is a critical value of $\Psi_K$. }

{\bf Proof.} For the reader's convenience, we sketch a brief proof here and refer to Sections V.2 and V.3
of \cite{Eke1} for related details.

By the proof of Theorem V.2.9 of \cite{Eke1}, if we replace $L_o^\beta$ and $\psi$ by $L^{2}({\bf R}/{\bf Z},{\bf R}^{2n})$
and $\Psi_K$ respectively, the Theorem V.2.9 of \cite{Eke1} also works.
Since the Fadell-Rabinowitz index  $\hat{I}$ has the properties of monotonicity, subadditivity, continuity which
are the only three properties of $I$ used in the proof of Proposition V.2.10 of \cite{Eke1}, then the
proof carries over verbatim of that of Proposition V.2.10 of \cite{Eke1}.\hfill\hb

Note that here we can't get $c_i\neq 0$ and prove Proposition 3.5 below, because it depends on Proposition 3.4 and
the identity (\ref{3.10}) below, but we should firstly prove Proposition 3.1 and 3.3 in order to get the identity (\ref{3.10})
by the method of Lemma V.3.8 of \cite{Eke1}.

{\bf Definition 3.2.} {\it Suppose $u$ is a nonzero critical
point of $\Psi_K$, and $\Nn$ is an $S^1$-invariant
open neighborhood of $S^1\cdot u$ such that
$crit(\Psi_K)\cap(\Lambda_K(u)\cap \Nn)=S^1\cdot u$. Then
the $S^1$-critical modules of $S^1\cdot u$ is defined by
\bea C_{S^1,\; q}(\Psi_K, \;S^1\cdot u)
=H_{q}((\Lambda_K(u)\cap\Nn)_{S^1},\;
((\Lambda_K(u)\setminus S^1\cdot u)\cap\Nn)_{S^1}),
\lb{3.3}\eea
where $\Lambda_K(u)=\{w\in L^{2}({\bf R}/{\bf Z},{\bf R}^{2n})\;|\;
\Psi_K(w)\le\Psi_K(u)\}$.}

Comparing with Theorem 4 in P.219 of \cite{Eke1}, we have the following

{\bf Proposition 3.3.} {\it For every $i\geq d(K)/2+1$, there exists a point
$u\in L^{2}({\bf R}/{\bf Z},{\bf R}^{2n})$ such that}
\bea
&& \Psi_K^\prime(u)=0,\quad \Psi_K(u)=c_i, \lb{3.4}\\
&& C_{S^1,\; 2(i-1)}(\Psi_K, \;S^1\cdot u)\neq 0. \lb{3.5}\eea

{\bf Proof.} By Lemma 8 in P.206 of \cite{Eke1}, we can use
Theorem 1.4.2 of \cite{Cha1} in the equivariant form  to obtain
\be H_{S^1,\,\ast}(\Psi_K^{c_i+\epsilon},\;\Psi_K^{c_i-\epsilon})
=\bigoplus_{\Psi_K(u)=c_i}C_{S^1,\; \ast}(\Psi_K, \;S^1\cdot u),\lb{3.6}\ee
for $\epsilon$  small enough such that the interval
$(c_i-\epsilon,\,c_i+\epsilon)$ contains no critical values of $\Psi_K$
except $c_i$.

Similar to P.431 of \cite{EkH1}, we have
\be H^{2(i-1)}((\Psi_K^{c_i+\epsilon})_{S^1},\,(\Psi_K^{c_i-\epsilon})_{S^1})
\mapright{q^\ast} H^{2(i-1)}((\Psi_K^{c_i+\epsilon})_{S^1} )
\mapright{p^\ast}H^{2(i-1)}((\Psi_K^{c_i-\epsilon})_{S^1}),\lb{3.7}\ee
where $p$ and $q$ are natural inclusions. Denote by
$f: (\Psi_K^{c_i+\epsilon})_{S^1}\rightarrow CP^\infty$ a classifying map and let
$f^{\pm}=f|_{(\Psi_K^{c_i\pm\epsilon})_{S^1}}$. Then clearly each
$f^{\pm}: (\Psi_K^{c_i\pm\epsilon})_{S^1}\rightarrow CP^\infty$ is a classifying
map on $(\Psi_K^{c_i\pm\epsilon})_{S^1}$. Let $\eta \in H^2(CP^\infty)$ be the first
universal Chern class.

By definition of $c_i$, we have $\hat I(\Psi_K^{c_i-\epsilon})< i$, hence
$(f^-)^\ast(\eta^{i-1})=0$. Note that
$p^\ast(f^+)^\ast(\eta^{i-1})=(f^-)^\ast(\eta^{i-1})$.
Hence the exactness of (\ref{3.7}) yields a
$\sigma\in H^{2(i-1)}((\Psi_K^{c_i+\epsilon})_{S^1},\,(\Psi_K^{c_i-\epsilon})_{S^1})$
such that $q^\ast(\sigma)=(f^+)^\ast(\eta^{i-1})$.
Since $\hat I(\Psi_K^{c_i+\epsilon})\ge i$, we have $(f^+)^\ast(\eta^{i-1})\neq 0$.
Hence $\sigma\neq 0$, and then
$$H^{2(i-1)}_{S^1}(\Psi_K^{c_i+\epsilon},\Psi_K^{c_i-\epsilon})=
H^{2(i-1)}((\Psi_K^{c_i+\epsilon})_{S^1},\,(\Psi_K^{c_i-\epsilon})_{S^1})\neq 0. $$
Thus the proposition follows from (\ref{3.6}) and the universal coefficient
theorem. \hfill\hb

Now we define two numbers $\gamma_\alpha^+(\Sigma)$ and $\gamma_\alpha^-(\Sigma)$ by:
\bea\gamma_\alpha^+(\Sigma)&=&\limsup_{i\rightarrow \infty}{[(-c_i)^{\frac{2-\alpha}{\alpha}}i]^{-1}},\nn\\
\gamma_\alpha^-(\Sigma)&=&\liminf_{i\rightarrow \infty}{[(-c_i)^{\frac{2-\alpha}{\alpha}}i]^{-1}},\nn\eea
and we set \bea \gamma^+(\Sigma)&=&\frac{\alpha}{4}(1-\frac{\alpha}{2})^{\frac{2-\alpha}{\alpha}}\gamma_\alpha^+(\Sigma),\lb{3.8}\\
\gamma^-(\Sigma)&=&\frac{\alpha}{4}(1-\frac{\alpha}{2})^{\frac{2-\alpha}{\alpha}}\gamma_\alpha^-(\Sigma).\lb{3.9}\eea
Then by the proofs of Lemma V.3.8 of \cite{Eke1} and Proposition 2.8 of \cite{LLW1}, noticing that when $\Sigma$ is convex, the Viterbo index
$i^v(y)$ and nullity $\nu^v(y)$ are the same as Ekeland index and nullity, we have
\bea \gamma^+(\Sigma_R)=\gamma^-(\Sigma_R)=\frac{\pi R^2}{2n},\lb{3.10}\eea
where $\Sigma_R$ is the sphere of radius $R$ in $\R^{2n}$.

{\bf Proposition 3.4.} {\it  We have $0<\gamma^-(\Sigma)\leq \gamma^+(\Sigma)<+\infty$ for any $\Sigma\in \H_{st}(2n)$.}

{\bf Proof.} Since $\Sigma$ is star-shaped, there exist some $0<r<R$ such that
\bea R^{-\alpha}|x|^\alpha\leq H_\alpha(x)\leq r^{-\alpha}|x|^\alpha.\lb{3.11}\eea
We denote the modified Hamiltonian functions of $R^{-\alpha}|x|^\alpha$ and
$r^{-\alpha}|x|^\alpha$ in Proposition 2.1 by $H_R(x)$ and $H_r(x)$ respectively, and we
can also choose the functions to satisfy\bea H_R(x)\leq H(x)\leq H_r(x), \forall x\in \R^{2n},\lb{3.12}\eea
where $H(x)$ is the modified Hamiltonian function of $H_\alpha(x)$.

Denote by $\Psi_K^r$ and $\Psi_K^R$ the corresponding dual action functionals defined in (\ref{2.7})
associated with the Hamiltonians $H_r$ and $H_R$ respectively,
then by (\ref{3.12}) we have \bea \Psi_K^r\leq\Psi_K\leq\Psi_K^R.\lb{3.13}\eea
Define \bea c_i(r)&=&\inf{\{\delta\mid \widehat{I}((\Psi_K^r)^{\delta-})\geq i\}}\nn\\
c_i(R)&=&\inf{\{\delta\mid \widehat{I}((\Psi_K^R)^{\delta-})\geq i\}}\nn\eea
Then $c_i(r)\leq c_i\leq c_i(R)$ from (\ref{3.13}). Thus by definitions (\ref{3.8}), (\ref{3.9}), we have
\bea &&\gamma^+(\Sigma_r)\leq \gamma^+(\Sigma)\leq \gamma^+(\Sigma_R),\lb{3.14}\\
&&\gamma^-(\Sigma_r)\leq \gamma^-(\Sigma)\leq \gamma^-(\Sigma_R).\lb{3.15}\eea
Hence by (\ref{3.10}), (\ref{3.14}) and (\ref{3.15}), we obtain
$\frac{\pi r^2}{2n}\leq\gamma^-(\Sigma)\leq \gamma^+(\Sigma)\leq \frac{\pi R^2}{2n}$. \hfill\hb

{\bf Proposition 3.5.} {\it If $c_i=c_j$ for some $d(K)/2+1\leq i<j$, then there are infinitely many geometrically
distinct closed characteristics on $\Sg$.}

{\bf Proof.} Note that by Proposition 3.4, we have $c_i\neq 0$, $i\geq d(K)/2+1$. Then by the same proof of
Proposition V.3.3 of \cite{Eke1}, we prove our proposition. \hfill\hb

Since every solution $(\tau,y)\in\T(\Sigma,\alpha)$ gives rise to a sequence
$\{z_m^y\}_{m\in\N}$ of solutions of the given period-1 problem (\ref{2.1}),
and a sequence $\{u_m^y\}_{m\in\N}$ of critical points of $\Psi_K$ defined by
\bea z_m^y(t)&=&(m\tau)^{-\frac{1}{2-\alpha}} y(m\tau t),\lb{3.16}\\
u_m^y(t)&=&-J(m\tau)^{-\frac{\alpha-1}{2-\alpha}} \dot{y}(m\tau t)+
K(m\tau)^{-\frac{1}{2-\alpha}} y(m\tau t)\lb{3.17}\eea
From the proof of Proposition 2.8 of \cite{LLW1}, we know that $\Psi_K(u_m^y)$ is independent
of $K$, together with (V.3.45) of \cite{Eke1}, it follows that
\bea \Psi_K(u_m^y)=-(1-\frac{\alpha}{2})\left(\frac{2}{\alpha}k A(y)\right)^{-\frac{\alpha}{2-\alpha}},\lb{3.18}\eea
where the action of a closed characteristic $(\tau,y)$ is
defined by \[
A(y)=\frac{1}{2}\int_0^{\tau}{(Jy\cdot \dot{y})}dt.\]

{\bf Corollary 3.6.} {\it We have $\lim\limits_{i\rightarrow\infty}{c_i}=0$ and for every
$i\in\N$, there exists $(\tau,y)\in\T(\Sigma,\alpha)$ and $m\in\N$ such that
\bea && \Psi_K^\prime(u_m^y)=0,\quad \Psi_K(u_m^y)=c_{i+d(K)/2},\lb{3.19} \\
&& i^v(y^m)\leq 2(i-1)\leq i^v(y^m)+\nu^v(y^m)-1,\lb{3.20}\eea
where $u_m^y$ is defined as in (\ref{3.17}). }

{\bf Definition 3.7.} {\it We call $(\tau,y)\in\T(\Sigma,\alpha)$ is $i$-essential
if there exists some $m\in\N$ such that (\ref{3.19}), (\ref{3.20}) hold. It is essential if it
is $i$-essential for some $i\in\N$. We denote by $\mathcal{C}$ the family
of essential closed characteristics on $\Sigma$.}

{\bf Theorem 3.8.} {\it We have $[1/\gamma^+(\Sigma), 1/\gamma^-(\Sigma)]\subseteq$
the closure of $\{\frac{\hat{i}(y)}{A(y)}\mid y\in \mathcal{C}\}$, where
$\hat{i}(y)\equiv \lim\limits_{m\rightarrow \infty}{\frac{i(y^m)}{m}}$ is the mean index of $(\tau,y)$. }

The proof of Theorem 3.8 relies on the following:

{\bf Lemma 3.9.} {\it There exists constant $d$,  which only depends on $\Sigma$, such that
whenever $y\in\mathcal{C}$ is $i$-essential for some $i\in\N$, we
have\bea \left|\frac{1}{C_\alpha}\frac{\hat{i}(y)}{A(y)}-
(d(K)/2+i)|c_{d(K)/2+i}|^{\frac{2-\alpha}{\alpha}}\right|\leq d|c_{d(K)/2+i}|^{\frac{2-\alpha}{\alpha}},\nn\eea
where $C_\alpha=\frac{4}{\alpha}(1-\frac{\alpha}{2})^{-\frac{2-\alpha}{\alpha}}$.}

{\bf Proof.} By Theorem 10.1.1, 10.1.2 of \cite{Lon2} and Lemma 2.4, the Viterbo index
$i^v(y)$ has the property of Proposition I.5.21 of \cite{Eke1}, note that Theorem V.1.4 of \cite{Eke1} also holds
for star-shaped hypersurfaces, then our lemma follows by the same proof of Lemma V.3.12
of \cite{Eke1}.\hfill\hb

{\bf Proof of Theorem 3.8.} From Lemma 3.9 instead of Lemma V.3.12
of \cite{Eke1}, our theorem follows by the same proof of Theorem V.3.11 of \cite{Eke1}.\hfill\hb

Now by the same proof of Theorem V.3.15 of \cite{Eke1}, we obtain

{\bf Theorem 3.10.} {\it If $\mathcal{C}$ is finite. Then we have
\bea &&\gamma(\Sigma)\equiv\gamma^+(\Sigma)=\gamma^-(\Sigma),\nn\\
&&\frac{\hat{i}(y)}{A(y)}=\frac{1}{\gamma(\Sigma)},\; \forall y\in\mathcal{C},\nn\\
&&\sum_{y\in\mathcal{C}}{\frac{1}{\hat{i}(y)}}\geq \frac{1}{2}.\lb{3.21}\eea}

By (\ref{3.21}), we have

{\bf Corollary 3.11.} {\it If there is a closed characteristic on $\Sg\in\H_{st}(2n)$ whose mean index is greater
than 2, then there exist at least two closed characteristics on $\Sg$.}

\subsection{ Critical values in the pinched case }
In this subsection, we prove under suitable pinching condition, the critical values $c_{i+d(K)/2}$ found
in Subsection 3.1 correspond to $n$ distinct closed characteristics for $1\leq i\leq n$.

Let $n(y)$ be the unit outward normal vector of $\Sigma$ at $y$ and
$d(y):=n(y)\cdot y$, i.e., the distance between the origin of $\R^{2n}$ and the tangent hyperplane to $\Sigma$
at $y$, then $d(y)>0$ for all $y\in \Sigma$ since $\Sigma$ is strictly star-shaped. Let $d=\min{\{d(y): y\in \Sigma\}}$,
$R=\max{\{|y|: y\in \Sigma\}}$. Then we have

\medskip

{\bf Theorem 3.12.} {\it Suppose that $\Sigma$ satisfies the pinching condition $R^2<2d^2$, then the critical values $c_{i+d(K)/2}$ found
in Proposition 3.1 and Corollary 3.6 correspond to at least $n$ distinct closed characteristics for $1\leq i\leq n$.}

\medskip

{\bf Proof.} We carry out our proof in two steps:

{\bf Step 1.} {\it We have\bea c_{i+d(K)/2}\leq -(1-\frac{\alpha}{2})(\frac{\alpha}{2\pi R^2})^{\frac{\alpha}{2-\alpha}}, \quad\forall\  1\leq i\leq n.\lb{3.22}\eea}

In fact, when $\Sigma=\Sigma_R$ is the sphere of radius $R$ in $\R^{2n}$,  by the proofs of Lemma V.3.8 of \cite{Eke1}
and Proposition 2.8 of \cite{LLW1}, noticing that when $\Sigma$ is convex, the Viterbo index
$i^v(y)$ and nullity $\nu^v(y)$ are the same as Ekeland index and nullity, we obtain that the corresponding critical values $c^R_{i+d(K)/2}$ found in
Proposition 3.1 and Corollary 3.6 satisfy\bea c^R_{i+d(K)/2}=-(1-\frac{\alpha}{2})(\frac{\alpha}{2\pi R^2})^{\frac{\alpha}{2-\alpha}},
\quad\forall 1\leq i\leq n,\nn\eea
which, together with (\ref{3.2}) and (\ref{3.13}), yields (\ref{3.22}).

{\bf Step 2.} {\it We have \bea c_{i+d(K)/2}\geq (1-\frac{\alpha}{2})(\frac{\alpha}{2\pi d^2})^{\frac{\alpha}{2-\alpha}},\quad \forall 1\leq i\leq n.\lb{3.23}\eea}

In fact, when we replace $r$ in the proof of Theorem V.1.4 of \cite{Eke1} by $d$, then
Theorem V.1.4 of \cite{Eke1} holds for star-shaped hypersurface $\Sigma$, i.e., for every closed characteristic $(\tau, y)$ on star-shaped hypersurface $\Sigma$,
there holds $A(\tau, y)\geq \pi d^2$, which, together with (\ref{3.18}), yields (\ref{3.23}).

Now, combining (\ref{3.22})-(\ref{3.23}), (\ref{3.18}) and $R^2<2d^2$, by Proposition 3.5 we obtain that the critical values $c_{i+d(K)/2}$
correspond to at least $n$ distinct closed characteristics for $1\leq i\leq n$.\hfill\hb

\setcounter{equation}{0}
\section{Proofs of Theorems 1.1 and 1.4-1.5}
In this section, we prove Theorems 1.1 and 1.4-1.5.

\medskip

{\bf Lemma 4.1.} (cf. Proposition 6.2 of \cite{LLo2}) {\it Let $\Sg\in \H_{st}(4)$ satisfy $\,^{\#}\T(\Sg)=2$. Denote
the two geometrically distinct prime closed characteristics by $\{(\tau_j,\; y_j)\}_{1\le j\le 2}$. If
$i(y_j)\geq 0, j=1, 2$, then both of the closed characteristics are elliptic.}

\medskip

{\bf Proof of Theorem 1.1.} Let $\Sg \in \H_{st}(4)$ satisfy $\,^{\#}\T(\Sg)=2$ and $R^2< 2d^2$, we denote by $\{(\tau_1,\; y_1), \;(\tau_2,\; y_2)\}$
the two geometrically distinct prime closed characteristics on $\Sg$, and by $\ga_j\equiv \ga_{y_j}$ the associated
symplectic paths of $(\tau_j,\; y_j)$ for $1\le j\le 2$. Then by Lemma 3.3 of \cite{HuL1} (cf. also Lemma 15.2.4 of
\cite{Lon2}), there exist
$P_j\in \Sp(4)$ and $M_j\in \Sp(2)$ such that
\be   \ga_j(\tau_j)=P_j^{-1}(N_1(1,\,1)\diamond M_j)P_j, \quad {\rm for}\;\;j=1, 2.  \lb{4.1}\ee
Note that by Section 9 of \cite{Vit2}, we know that there exists at least one non-hyperbolic closed characteristic
on $\Sg$ and it is certainly elliptic when $n=2$. In the following, we prove Theorem 1.1 by contradiction. Without
loss of generality, we assume that $(\tau_1, y_1)$ is elliptic and $(\tau_2, y_2)$ is hyperbolic.

For these two closed characteristics, we have the following properties:

\medskip

{\bf Claim 1.} {\it The closed characteristics $(\tau_1, y_1)$, $(\tau_2, y_2)$ satisfy

(i) $i(y_2^m) = m(i(y_2)+3) - 3$ and $\nu(y_2^m)=1$, $\forall$ $m\in\N$, and thus $\hat{i}(y_2) = i(y_2) + 3$.

(ii) $\hat{i}(y_2) >0$.

(iii) $\hat{i}(y_1) \in \Q$.

(iv) If $i(y_2)$ is even, then $i(y_2^2)-i(y_2)\in 2\Z-1$, $\hat{\chi}(y_2)=\frac{1}{2}$, and $i(y_2)\geq -2$.}

\medskip

In fact, by Theorem 8.3.1 of \cite{Lon2}, we have $i(y_2,m) = m(i(y_2, 1)+1)-1$, $\forall\ m\in\N$. Together with Lemma 2.4, we obtain (i).

We claim $\hat{i}(y_2)\neq 0$. In fact, because $y_2$ is hyperbolic, $y_2^m$ is non-degenerate for every
$m\ge 1$. Thus if $\hat{i}(y_2)=0$, we then have $i(y_2^m)=i(y_2,m)-2=-3$ for all $m\ge 1$. Then the Morse-type
number satisfies $m_{-3}=+\infty$. But then $\hat{i}(y_1)$ must be positive by Theorem 5.6, and
contributions of $\{y_1^m\}$ to every Morse-type number thus must be finite. Then the Morse inequality yields a
contradiction and proves the claim (cf. the proof below (9.3) of \cite{Vit2} for details).

If $\hat{i}(y_2)<0$, by (\ref{5.24}) we obtain
\be    \frac{\hat{\chi}(y_2)}{\hat{i}(y_2)} = 0.   \lb{4.2}\ee
But because $(\tau_2, y_2)$ is hyperbolic, by (\ref{5.22}) we have $\hat{\chi}(y_2)\not= 0$, which
contradicts to (\ref{4.2}) and proves (ii).

If $(\tau_1, y_1)$ and its iterates are all non-degenerate, since $(\tau_1, y_1)$ is elliptic, then
$\hat{i}(y_1)$ must be irrational by Corollary 8.3.2 of \cite{Lon2} and then so is
$\frac{\hat{\chi}(y_1)}{\hat{i}(y_1)}$, because $\hat{\chi}(y_1)$ is rational and nonzero by (\ref{5.22}). Then by (\ref{5.23})
of Theorem 5.6, the other closed characteristic $(\tau_2, y_2)$ must possess an irrational mean index
$\hat{i}(y_2)$, which contradicts to the second identity in (i), and thus $\hat{i}(y_1)$ must be
rational, which proves (iii).

If $i(y_2)$ is even, then $i(y_2^2)-i(y_2)=i(y_2)+3\in 2\Z-1$ by (i), which together (\ref{5.22}) implies $\hat{\chi}(y_2)=\frac{1}{2}$.
Then $i(y_2)\geq -2$ follows from (i) and (ii).

The proof of Claim 1 is complete.

\medskip

By (iii) of Claim 1, we only need to consider the following four cases according to the classification
of basic norm forms of $\ga_1(\tau_1)$. In the following we use the notations from Definition 1.8.5
and Theorem 1.8.10 of \cite{Lon2}, and specially we let $R(\th) =  \left(\begin{array}{cc}
                                                             \cos{\th} & -\sin{\th} \\
                                                             \sin{\th} &  \cos{\th} \\
                                                                \end{array}\right)$ with $\th\in\R$,
and use $M\dm N$ to denote the symplectic direct sum of two symplectic matrices $M$ and $N$ as
in pages 16-17 of \cite{Lon2}.

\medskip

{\bf Case 1.} {\it $\ga_1(\tau_1)$ can be connected to $N_1(1,1)\dm N_1(-1,b)$ within
$\Om^0(\ga_1(\tau_1))$ with $b=0$ or $\pm 1$.}

\medskip

In this case, by Theorems 8.1.4 and 8.1.5 of \cite{Lon2}, and Lemma 2.4, we have
\be  i(y_1, 1) \quad {\rm and}\quad i(y_1) \quad {\rm are\;even}.  \lb{4.3}\ee
By Theorem 1.3 of \cite{Lon1}, we have
\bea
i(y_1,m) &=& m(i(y_1, 1)+1)-1, \quad {\rm for}\;\; b=1;\nn\\
i(y_1, m) &=& m(i(y_1,1)+1)-1-\frac{1+(-1)^m}{2}, \quad {\rm for}\;\;b=0, -1. \nn\eea
By Lemma 2.4, we obtain
\bea
i(y_1^m) &=& m(i(y_1)+3)-3, \quad {\rm for}\;\; b=1;   \lb{4.4}\\
i(y_1^m) &=& m(i(y_1)+3)-3-\frac{1+(-1)^m}{2}, \quad {\rm for}\;\;b=0, -1. \lb{4.5}\eea
Then in both cases we obtain
\be \hat{i}(y_1) = i(y_1)+3.  \lb{4.6}\ee

Next we separate our proof in two subcases according to the parity of $i(y_2)$.

\medskip

{\bf Subcase 1.1.} {\it $i(y_2)$ is even.}

\medskip

By (iv) of Claim 1, we have $i(y_2)\geq  0$ or $i(y_2)=-2$.
We continue our proof in two steps according to the value of $i(y_2)$:

\medskip

{\bf Step 1.1.} $i(y_2)\geq  0$.

\medskip

In this step,  by  (i) of Claim 1 we have \bea \hat{i}(y_2)\geq 3, \lb{4.7}\eea
which together with (iv) of Claim 1 implies
\bea \frac{\hat{\chi}(y_2)}{\hat{i}(y_2)}\leq\frac{1}{6}.\lb{4.8}\eea
Combining (\ref{4.8}) with Theorem 5.6, we obtain \bea\hat{i}(y_1)> 0,\quad \frac{\hat{\chi}(y_1)}{\hat{i}(y_1)}=\frac{1}{2}
-\frac{\hat{\chi}(y_2)}{\hat{i}(y_2)}\geq \frac{1}{3}. \lb{4.9}\eea

Note that by Proposition 5.4 and the form of $\ga_1(\tau_1)$, we have $K(y_1)=2$. Thus by (\ref{5.21})  and (\ref{4.9}), we obtain
\be   0 < \hat{\chi}(y_1)=\frac{1+(-1)^{i(y_1^2)}(k_0(y_1^2)-k_1(y_1^2)+k_2(y_1^2))}{2}.   \lb{4.10}\ee
Since at most one of $k_l(y_1^2)$s for $0\leq l\leq 2$ can be non-zero by (iv) of Remark 5.5, we obtain
\be   (-1)^{i(y_1^2)+l}k_l(y_1^2)\ge 0, \qquad {\rm for}\;\;l=0, 1, 2.   \lb{4.11}\ee

When $m$ is odd, we have $\nu(x_1^m)=1$ by the assumption on $\ga_1(\tau_1)$. In this case, because
$i(y_1)$ is even by (\ref{4.3}), we have $i^v(x_1^m) = i(y_1^m) = m(i(y_1)+3)-3$ is even, and then
$$  \bb(x^m) = (-1)^{i^v(x_1^m)-i^v(x_1)} = 1, $$
where we denote by $x_j$ the critical point of $F_{a,K}$ corresponding to $y_j$ for $j=1$ and $2$.
Thus by (\ref{5.17}) of Proposition 5.2 for every odd $m\in\N$, we obtain
\bea
C_{S^1,\;d(K)+k}(F_{a,K},\;S^1\cdot {x}_1^{m})=\Q, && {\rm if}\;\;k=i(y_1^m),   \lb{4.12}\\
C_{S^1,\;d(K)+k}(F_{a,K},\;S^1\cdot {x}_1^{m})=0, && {\rm if}\;k\not=i(y_1^m),   \lb{4.13}\eea
where (\ref{4.13}) holds specially when $k\in 2\Z-1$.

When $m$ is even, we consider two cases: (A-1) for $b=0, -1$ with (\ref{4.5}); (B-1) $b=1$ with (\ref{4.4}).

\medskip

{\bf (A-1)} {\it $m$ is even, $b=0$ or $-1$, and (\ref{4.5}) holds.}

\medskip

In this case, $i(y_1^2)$ is even by (\ref{4.5}). Therefore by (\ref{4.10})-(\ref{4.11}) we obtain
\be  k_1(y^2) = 0, \quad \hat{\chi}(y_1)=\frac{1+(k_0(y_1^2)+k_2(y_1^2))}{2}>0.  \lb{4.14}\ee
Because $K(y_1)=2$, we then obtain
\be C_{S^1,\;d(K)+2k-1}(F_{a,K},\;S^1\cdot {x}_1^{m})=0, \qquad\forall k\in\Z, m\in 2\N.   \lb{4.15}\ee
Therefore, when $b=0, -1$, from (\ref{4.12}), (\ref{4.13}) and (\ref{4.15}) we obtain
\be C_{S^1,\;d(K)+2k-1}(F_{a,K},\;S^1\cdot {x}_1^{m}) = 0, \qquad \forall k\in\Z, m\in \N. \lb{4.16}\ee

\medskip

{\bf (B-1)} {\it $m$ is even, $b=1$, and (\ref{4.4}) holds.}

\medskip

In this case, $i(y_1^2)$ is odd by (\ref{4.4}). Therefore by (\ref{4.10})-(\ref{4.11}) we obtain
\be  k_0(y^2) = k_2(y^2) = 0, \quad 0<\hat{\chi}(y_1)=\frac{1+k_1(y_1^2)}{2}.  \lb{4.17}\ee
Because $K(y_1)=2$, we then obtain
\be C_{S^1,\;d(K)+2k-1}(F_{a,K},\;S^1\cdot {x}_1^{m})=0, \qquad \forall k\in\Z, m\in 2\N.   \lb{4.18}\ee
Therefore when $b=1$, from (\ref{4.12}), (\ref{4.13}) and (\ref{4.18}), we obtain
\be C_{S^1,\;d(K)+2k-1}(F_{a,K},\;S^1\cdot {x}_1^{m}) = 0, \qquad \forall k\in\Z, m\in \N. \lb{4.19}\ee

In summary, from (\ref{4.16}) and (\ref{4.19}), for any case we have
\be C_{S^1,\;d(K)+2k-1}(F_{a,K},\;S^1\cdot {x}_1^{m}) = 0, \qquad \forall k\in\Z, m\in \N. \lb{4.20}\ee

Note that in Subcase 1.1, $i(y_2)$ is even and $(\tau_2, y_2)$ is hyperbolic, then by
(\ref{5.17}) of Proposition 5.2, we obtain
\be C_{S^1,\;d(K)+2k-1}(F_{a,K},\;S^1\cdot {x}_2^{m}) = 0, \qquad \forall k\in\Z, m\in \N. \lb{4.21}\ee
Combining (\ref{4.20}) and (\ref{4.21}), we have $m_{2q-1}=0$ for every $q\in\Z$ and $U(t)\equiv 0$ in (\ref{5.28}).
Here and below in this Section $m_i$ denotes the coefficient of $t^i$ of $M(t)=\sum_{i\in \Z}{m_it^i}$ in
(\ref{5.28}). Then \bea  \sum_{i\in\Z}{m_it^i}=\frac{1}{1-t^2}=\sum_{i\in\N}{t^{2i-2}}.  \nn\eea
Thus $i(y_j)\geq 0$ for $j=1,2$ by Proposition 5.2. By Lemma 4.1 we know that the two closed characteristics are elliptic, which contradicts to our assumption.

\medskip

{\bf Step 1.2.} $i(y_2)=-2$.

\medskip

In this step,  by (i) of Claim 1 we have \bea \hat{i}(y_2)=1,\lb{4.22}\eea
which together with (iv) of Claim 1 implies \bea\frac{\hat{\chi}(y_2)}{\hat{i}(y_2)}=\frac{1}{2}.\lb{4.23}\eea

We continue in two cases: (A-2) $i(y_1)=-2$; (B-2) $i(y_1)\leq-4$ or $i(y_1)\geq 0$.

\medskip

{\bf (A-2)} {\it $i(y_1)=-2$ holds.}

\medskip

In this case, by (\ref{4.6}) we have\bea\hat{i}(y_1)=1,\lb{4.24}\eea
Combining (\ref{4.22}), (\ref{4.24}) with Theorem 5.6, we obtain \bea\frac{\hat{\chi}(y_1)}{\hat{i}(y_1)}
+\frac{\hat{\chi}(y_2)}{\hat{i}(y_2)}=\frac{1}{2},\lb{4.25}\eea
which together with (\ref{4.24}) implies \bea\hat{\chi}(y_1)=0.\lb{4.26}\eea

Since $i(y_1)=-2$ and $i(y_2)=-2$ in this case, there hold $i(c_k^m)\ge -2$, $\forall\ m\ge 1$, $k=1,2$. Then by
(\ref{5.16})-(\ref{5.17}) of Proposition 5.2, (\ref{4.4})-(\ref{4.5}) and (i) of Claim 1, we obtain
$m_{-2}\geq 2$ and $m_q=0$ for $q<-2$, which together with (\ref{5.28}) gives \bea m_{-1}=u_{-2}+u_{-1}=m_{-2}+u_{-1}\geq 2. \lb{4.27}\eea

On the other hand, by (i) of Claim 1,
$i(y_2^m) = m-3$, then by (\ref{5.17}) of Proposition 5.2, we get
\bea C_{S^1,\;d(K)-1}(F_{a,K},\;S^1\cdot {x}_2^{m})=0,\quad \forall
m\in\N,\lb{4.28}\eea
which, together with (\ref{4.27}), implies that $y_1^m$ has contribution to $m_{-1}$ for some $m\in\N$
and $y_2^m$ has no contribution to $m_{-1}$ for all $m\in\N$. Note that by (\ref{4.4})-(\ref{4.5}) and
Proposition 5.2, then we have $m_{-1}=k_0(y_1^2)\ge 2$ when $b=1$ and (\ref{4.4}) holds, or
$m_{-1}=k_1(y_1^2)\ge 2$ when $b=0$ or $-1$, and (\ref{4.5}) holds. Together with (iv) of Remark 5.5, (\ref{4.27}) and (\ref{5.21}), it yields
\bea \hat{\chi}(y_1)=\frac{1+(-1)^{i(y_1^2)}(k_0(y_1^2)-k_1(y_1^2)+k_2(y_1^2))}{2}\neq 0,\nn\eea
which contradicts to (\ref{4.26}).

\medskip

{\bf (B-2)} {\it $i(y_1)\leq-4$ or $i(y_1)\geq 0$ holds.}

\medskip

In this case, by (\ref{4.4})-(\ref{4.5}) and
Proposition 5.2, $y_1^m$ has no contribution to $m_{-1}$ and $m_{-3}$ for all $m\in\N$.
Note that $i(y_2^m) = m-3$ in Step 2, then by Proposition 5.2, $y_2^m$ has no contribution to $m_{-1}$ and $m_{-3}$ for all $m\in\N$ too.
Thus $m_{-1}=m_{-3}=0$, which together with (\ref{5.28}) gives \bea m_{-2}=0. \lb{4.29}\eea
Since $i(y_2)=-2$, then by Proposition 5.2, we have $C_{S^1,\;d(K)-2}(F_{a,K},\;S^1\cdot {x}_2)=\Q$ which implies
$m_{-2}>0$, it contradicts to (\ref{4.29}).

\medskip

{\bf Subcase 1.2.} {\it $i(y_2)$ is odd.}

\medskip

In this subcase, by (i)-(ii) of Claim 1, we have $i(y_2)\geq -1$ and it is odd.
When $i(y_2)\geq 1$, then by the same proof of Subcase 1.1 of Theorem 1.4 in \cite{LLo2}, we get a contradiction.
Thus we can assume that $i(y_2)=-1$. Then by (i) of Claim 1, we have
$i(y_2^m) = 2m-3$, it together with Proposition 5.2 gives $C_{S^1,\;d(K)+2i-2}(F_{a, K},\;S^1\cdot {x}_2^m)=0$, $\forall\ m\ge 1$, $i=1,2$.
On the other hand, by Proposition 3.3 and Theorem 3.12, we have $C_{S^1,\;d(K)+2i-2}(\Psi_K,\;S^1\cdot {u}^{y_2})\neq0$ for $i=1$ or $i=2$,
where ${u}^{y_2}$ is the critical point of $\Psi_K$ corresponding to $y_2$. By the same proof of Proposition 3.6 of \cite{Wan1},
we have $C_{S^1,\;d(K)+2i-2}(\Psi_K,\;S^1\cdot {u}^{y_2})\cong C_{S^1,\;d(K)+2i-2}(F_{a, K},\;S^1\cdot {x}_2)$, it is a contradiction.

\medskip

{\bf Case 2.} {\it $\ga_1(\tau_1)$ can be connected to $N_1(1,1)\dm R(\th)$ within $\Om^0(\ga_1(\tau_1))$
with some $\th\in (0,\pi)\cup (\pi,2\pi)$ and $\th/\pi\in \Q$.}

\medskip

In this case, we have always $K(y_1)\ge 3$ by the definition of $\th$. By Theorems 8.1.4 and 8.1.7
of \cite{Lon2} and Lemma 2.4 we obtain
\be  i(y_1,1)\quad {\rm and}\quad i(y_1) \quad {\rm are\;even}.  \lb{4.30}\ee

By Theorem 1.3 of \cite{Lon1} (i.e., Theorem 8.3.1 of \cite{Lon2}), we have
$$   i(y_1,m) = mi(y_1,1) + 2E(\frac{m\th}{2\pi}) - 2, \quad\forall\ m\ge 1, $$
which, together with Lemma 2.4, yields
\be  i(y_1^m) = m(i(y_1)+2) + 2E(\frac{m\th}{2\pi}) - 4, \quad\forall\ m\ge 1.   \lb{4.31}\ee
Then
\be   \hat{i}(y_1) = i(y_1) + 2 + \frac{\th}{\pi}.   \lb{4.32}\ee

We have two subcases according to the parity of $i(y_2)$.

\medskip

{\bf Subcase 2.1.} {\it $i(y_2)$ is odd.}

\medskip

For this case, as the same proof of Case 2 of Theorem 1.4 in \cite{LLo2}, we can get a contradiction.

\medskip

{\bf Subcase 2.2.} {\it $i(y_2)$ is even.}

\medskip

Next we continue our proof in two steps according to the value of $i(y_1)\in 2\Z$ by (\ref{4.30}).

\medskip

{\bf Step 2.1.} $i(y_1)\leq-4$.

\medskip

In this step, by (\ref{4.32}), we have $\hat{i}(y_1)<0$, together with Theorem 5.6, we obtain
\bea \frac{\hat{\chi}(y_2)}{\hat{i}(y_2)}=\frac{1}{2},\lb{4.33}\eea
and \bea \frac{\hat{\chi}(y_1)}{\hat{i}(y_1)}=0, \lb{4.34}\eea
which implies \bea \hat{\chi}(y_1)=0. \lb{4.35}\eea

By (\ref{4.33}) and (iv) of Claim 1, we get $\hat{i}(y_2)=1$, which, together with (i) of Claim 1, implies
\bea i(y_2)=-2. \lb{4.36}\eea
By (\ref{5.21}), we have \be  \hat{\chi}(y_1)
\;=\; \frac{K(y_1)-1+k_0(y_1^{K(y_1)})-k_1(y_1^{K(y_1)})+k_2(y_1^{K(y_1)})}{K(y_1)},
             \lb{4.37}\ee
which, together with (iv) of Remark 5.5 and (\ref{4.35}), yields
\bea k_1(y_1^{K(y_1)})=K(y_1)-1>0,\quad k_0(y_1^{K(y_1)})=k_2(y_1^{K(y_1)})=0. \lb{4.38}\eea
Since $i(y_1^{mK(y_1)})\le(-2+\frac{\theta}{\pi})mK(y_1)-4<-4$, $\forall m\in\N$, then
by Proposition 5.2 and (\ref{4.38}) we know that $y_1^{mK(y_1)}$ has no contribution to $m_{-1}$ and $m_{-3}$. On the other hand,
note that $i(y_1^{m})$ is even, it follows from Proposition 5.2 that $y_1^{m}$ has no contribution to $m_{-1}$ and $m_{-3}$ for $m\not= 0\  ({\rm mod}\; K(y_1))$.
In addition, $y_2^{m}$ also has no contribution to $m_{-1}$ and $m_{-3}$ since $y_2$ is hyperbolic and $i(y_2)\in2\Z$. Hence we obtain $m_{-1}=m_{-3}=0$, which, together with (\ref{5.28}), yields $m_{-2}=0$. But by (\ref{4.36}) and Proposition 5.2, $y_2$ contributes $1$ to $m_{-2}$. So we get a contradiction.

\medskip

{\bf Step 2.2.} $i(y_1)\geq-2$.

\medskip

In this step, note that $i(y_1)$ and $i(y_2)$ are even, we have either $i(y_j)\geq 0$ for $1\leq j\leq 2$, or
$i(y_j)=-2$ for some $1\leq j\leq 2$, if the former holds, then by Lemma 4.1, $y_1$ and $y_2$ are elliptic which
contradicts to our assumption. Thus we can assume that $i(y_j)=-2$ for some $1\leq j\leq 2$. Then by Proposition 5.2
we have \bea m_{-2}\geq1.\lb{4.39}\eea
Note that $i(y_j)\geq-2$ for $1\leq j\leq 2$ by (iv) of Claim 1. By (i) of Claim 1 and (\ref{4.31}) we have
\bea i(y_j^m)\geq-2,\quad\forall\ m\ge 1,\quad j=1,2.\lb{4.40}\eea
Thus we have
$$  C_{S^1,\;d(K)-q}(F_{a,K},\;S^1\cdot {x}_j^{m})=0, \qquad \forall\;m\in\N, \;q\ge 3,\; j=1,2.  $$
Hence we have
\bea  m_{-q}=0, \qquad \forall\;q\ge 3,  \nn\eea
which, together with (\ref{4.39}) and (\ref{5.28}), yields
\be   m_{-1}=u_{-1}+u_{-2}\ge u_{-2}=u_{-2}+u_{-3}= m_{-2}\ge 1.   \lb{4.41}\ee

Note that $y_2^{m}$ has no contribution to $m_{2k-1}$ for $k\in\Z$ since $y_2$ is hyperbolic and $i(y_2)\in2\Z$,
so some $y_1^m$ must have contribution to $m_{-1}$. Also note that, in this case, $C_{S^1,\;d(K)+2i-1}(F_{a, K},\;S^1\cdot {x}_1^m)=0$ for any $i\in\Z$ and $m\neq0\;(\mod\ K(y_1))$.
Therefore $m_{-1}$ only can be contributed by iterates $y_1^{K(y_1)}$. This implies $i(y_1^m)\le i(y_1^{K(y_1)})= -2$, $\forall\ 1\leq m\leq K(y_1)-1$.
Thus by (\ref{4.40}) and Proposition 5.2, we have
\bea i(y_1^m)=-2,\quad \forall\  1\leq m\leq K(y_1)-1,\quad m_{-1}=k_1(y_1^{K(y_1)}).\lb{4.42}\eea
By (\ref{4.41})-(\ref{4.42}), we have $k_1(y_1^{K(y_1)})=m_{-1}\geq m_{-2}\geq K(y_1)-1$, which together
with (\ref{4.37}) yields\bea \hat{\chi}(y_1)\leq 0. \lb{4.43}\eea

Noticing that $i(y_j)\geq-2$ for $1\leq j\leq 2$, then by (ii) of Claim 1 and (\ref{4.32}) we have $\hat{i}(y_j)>0$,
which, together with (\ref{4.43}) and Theorem 5.6, yields \bea
\frac{\hat{\chi}(y_2)}{\hat{i}(y_2)}=\frac{1}{2}-\frac{\hat{\chi}(y_1)}{\hat{i}(y_1)}\geq \frac{1}{2}.\lb{4.44}\eea
On the other hand, by (i) and (iv) of Claim 1, we have \bea
\frac{\hat{\chi}(y_2)}{\hat{i}(y_2)}=\frac{1}{2(i(y_2)+3)}\leq \frac{1}{2},\nn\eea
which together with (\ref{4.44}) implies \bea \hat{\chi}(y_1)=0,\qquad i(y_2)=-2. \lb{4.45}\eea
Then $y_2^m$ contributes exactly 1 to $m_{-2}$, which, together with (\ref{4.42}), yields $m_{-2}=K(y_1)$. Thus
by (\ref{4.37}), (\ref{4.41})-(\ref{4.42}), we obtain $\hat{\chi}(y_1)< 0$, which contradicts to (\ref{4.45}).

\medskip

{\bf Case 3.} {\it $\ga_1(\tau_1)$ can be connected to $N_1(1,1)\dm N_1(1,b)$ within $\Om^0(\ga_1(\tau_1))$
with $b=0$ or $1$.}

\medskip

In this case, we have $K(y_1)=1$ by Proposition 5.4, $i(y_1,1)$ and then $i(y_1)$ is even by Theorem 8.1.4 of \cite{Lon2} and Lemma 2.4. By
Theorem 8.3.1 of \cite{Lon2}, we obtain $i(y_1,m)=m(i(y_1, 1)+2)-2$ for all $m\in\N$. Thus by Lemma 2.4
we have
\bea i(y_1^m)=m(i(y_1)+4)-4, \quad \forall\;m\in\N,\qquad \hat{i}(y_1)=i(y_1)+4.  \lb{4.46}\eea
Then we can assume that $\hat{i}(y_1)\neq 0$, and $\{y_1^m\}_{m\in\N}$ has contributions
to the Morse-type numbers $\{m_q\}_{q\in\Z}$, which implies that exactly one of $k_l(y_1^m)$ for $0\leq l\leq 2$ is nonzero by (iv) of Remark 5.5.

In fact, if $y_1^m$ has no contribution
to any Morse-type number $m_{q}$, by the proof of Theorem 1.1 of \cite{LLo1} we obtain three closed characteristics, which
contradicts to our assumption. If $\hat{i}(y_1)= 0$ and $\{y_1^m\}_{m\in\N}$ has contributions
to the Morse type numbers $\{m_q\}_{q\in\Z}$, then $i(y_1^m)=-4$ by (\ref{4.46}) and
exactly one of $k_l(y_1^m)$ for $0\leq l\leq 2$ is nonzero by (iv) of Remark 5.5, then as the proof of (i) of Claim 1, we can get a contradiction.

Next we consider two subcases according to the parity of $i(y_2)$.

\medskip

{\bf Subcase 3.1.} {\it $i(y_2)$ is odd.}

\medskip

For this case, as the same proof of Case 3 of Theorem 1.4 in \cite{LLo2}, we can get a contradiction.

\medskip

{\bf Subcase 3.2.} {\it $i(y_2)$ is even.}

\medskip

In this case, we note that $\hat{i}(y_1)\neq0$. If $\hat{i}(y_1)<0$, by Theorem 5.6, we have \bea\frac{\hat{\chi}(y_1)}{\hat{i}(y_1)}=0.\nn\eea
Then $\hat{\chi}(y_1)= 0$. But exactly one of $k_l(y_1^m)$ for $0\leq l\leq 2$ is nonzero, which together with (\ref{5.21}) implies
\bea \hat{\chi}(y_1)=k_0(y_1)-k_1(y_1)+k_2(y_1)\neq 0,\nn\eea
which is a contradiction. Thus \bea \hat{i}(y_1)>0.\lb{4.47}\eea
Then by Theorem 5.6 and (i) and (iv) of Claim 1, we have \bea\frac{\hat{\chi}(y_1)}{\hat{i}(y_1)}=
\frac{1}{2}-\frac{\hat{\chi}(y_2)}{\hat{i}(y_2)}\geq 0.\nn\eea
Then there holds $\hat{\chi}(y_1)\geq 0$ by (\ref{4.47}). By (\ref{5.21}), we get
\bea k_0(y_1)-k_1(y_1)+k_2(y_1)=\hat{\chi}(y_1)\geq 0,\nn\eea
which, together with (iv) of Remark 5.5, yields $k_1(y_1)=0$. Then by Proposition 5.4 we know that $y_1^m$ has no contribution
to $m_{2q-1}$ for all $q\in\Z$. Note that $i(y_2)$ is even and $y_2$ is hyperbolic, then by Proposition 5.2,
$y_2^m$ also has no contribution to $m_{2q-1}$ for all $q\in\Z$. Thus $m_{2q-1}=0$ for every $q\in\Z$, which implies $U(t)\equiv 0$ in (\ref{5.28}),
then \bea  \sum_{i\in\Z}{m_it^i}=\frac{1}{1-t^2}=\sum_{i\in\N}{t^{2i-2}}.  \lb{4.48}\eea
Thus $i(y_2)\geq 0$ by Proposition 5.2 and (\ref{4.48}). Note that $i(y_1)\geq0$ or $i(y_1)=-2$ by (\ref{4.46})-(\ref{4.47}) and the fact that $i(y_1)$ is even. If
$i(y_1)\geq0$, then by Lemma 4.1 we know that the two closed characteristics are elliptic, which contradicts to our assumption.
Thus we suppose $i(y_1)=-2$. By (i) of Claim 1, it is impossible that $y_2^m$ contributes $1$ to every Morse-type number $m_q$ for $q\in2\N_0$. Noticing that exactly one of $k_l(y_1^m)$ for $0\leq l\leq 2$ is nonzero, by (\ref{4.48})
we have $k_2(y_1^m)=1$, which implies that $\{y_1^m\}_{m\in\N}$ contributes exactly 1
to every Morse-type number $m_q$ for $q\in2\N_0$, but $y_2$ also has contribution to some Morse-type number $m_q$, which contradicts to (\ref{4.48}).

\medskip

{\bf Case 4.} {\it $\gamma_1(\tau_1)$ can be connected to
$\left(\begin{array}{cc}
      1 & 1 \\
       0 & 1 \\
       \end{array} \right)\diamond\left(\begin{array}{cc}
                                          1 & -1 \\
                                            0 & 1 \\ \end{array}\right)$ within $\Omega^0(\gamma_1(\tau_1))$.}
\medskip

In this case, we have $K(y_1)=1$ by Proposition 5.4, $i(y_1,1)$ and then $i(y_1)$ is odd by Theorem 8.1.4 of \cite{Lon2} and Lemma 2.4. By
Theorem 8.3.1 of \cite{Lon1}, we have $i(y,m)=m(i(y, 1)+1)-1$ for all $m\in\N$. Thus by Lemma 2.4, we
obtain
\bea i(y_1^m)=m(i(y_1)+3)-3, \quad \forall\;m\in\N,\qquad \hat{i}(y_1)=i(y_1)+3. \lb{4.49}\eea
Then as Case 3, we can suppose that $\hat{i}(y_1)\neq 0$ and exactly one of $k_l(y_1^m)$ for $0\leq l\leq 2$ is nonzero.

We have two subcases according to the parity of $i(y_2)$.

\medskip

{\bf Subcase 4.1.} {\it $i(y_2)$ is odd.}

\medskip

For this case, as the same proof of Case 4 of Theorem 1.4 in \cite{LLo2} we can get a contradiction.

\medskip

{\bf Subcase 4.2.} {\it $i(y_2)$ is even.}

\medskip

As Subcase 3.2 we have \bea \hat{i}(y_1)>0,\qquad \hat{\chi}(y_1)\geq 0.\lb{4.50}\eea Then by (\ref{5.21}), we have
\bea -k_0(y_1)+k_1(y_1)=\hat{\chi}(y_1)\geq 0,\nn\eea
which, together with (iv) of Remark 5.5, implies $k_0(y_1)=0$. Then by Proposition 5.4 we know that $y_1^m$ has no contribution
to $m_{2q-1}$ for all $q\in\Z$. Note that $i(y_2)$ is even and $y_2$ is hyperbolic, then by Proposition 5.2,
$y_2^m$ also has no contribution to $m_{2q-1}$ for all $q\in\Z$. Thus $m_{2q-1}=0$ for every $q\in\Z$, which implies $U(t)\equiv 0$ in (\ref{5.28}),
then we have \bea  \sum_{i\in\Z}{m_it^i}=\frac{1}{1-t^2}=\sum_{i\in\N}{t^{2i-2}}.  \lb{4.51}\eea
Thus $i(y_2)\geq 0$ by Proposition 5.2 and (\ref{4.51}). Note that $i(y_1)\geq 1$ or $i(y_1)=-1$ by (\ref{4.49})-(\ref{4.50}) and the fact that $i(y_1)$ is odd. If
$i(y_1)\geq1$, then by Lemma 4.1 we know that the two closed characteristics are elliptic, which contradicts to our assumption.
Thus we suppose $i(y_1)=-1$. By (i) of Claim 1, it is impossible that $y_2^m$ contributes $1$ to every Morse-type number $m_q$ for $q\in2\N_0$. Noticing that exactly one of $k_l(y_1^m)$ for $0\leq l\leq 1$ is nonzero, by (\ref{4.51})
we have $k_1(y_1^m)=1$ and $\{y_1^m\}_{m\in\N}$ contributes exactly 1
to every Morse-type number $m_q$ for $q\in2\N_0$, but $y_2$ also has contribution to some Morse-type numbers $m_q$, which contradicts to (\ref{4.51}).

The proof of Theorem 1.1 is complete. \hfill\hb

\medskip

In the following, we explain why Theorems 1.4 and 1.5 hold.

\medskip

{\bf Proof of Theorem 1.4.} By Lemma 2.4 and Corollary 3.6, we have a similar result
as Theorem 1.6 of \cite{LoZ1}, and by Theorem 3.10, Lemma 3.1 of \cite{LoZ1} holds for
star-shaped hypersurfaces. Note that all the proofs of \cite{LoZ1} are based on the fact that
every closed characteristic $(\tau,y)$ on the hypersurface $\Sg$ in $\R^{2n}$ satisfies
$i(y)\geq n$ and Theorem 1.6, Lemma 3.1 of \cite{LoZ1}. Hence for dynamically convex star-shaped case,
all the theories of \cite{LoZ1} hold. Then combining it with Theorem 1.1 of \cite{Wan2} and Theorem 1.1 of \cite{HuO}, we get the desired results.\hfill\hb

\medskip

{\bf Proof of Theorem 1.5.} Note that all the proofs of \cite{WHL1} and \cite{Wan1} rely on
the resonance identity in Theorem 1.2 of \cite{WHL1}, the periodic property of critical
modules in Proposition 3.13 of \cite{WHL1}, and the results in \cite{LoZ1}.
Since we have extended the theories of \cite{WHL1} to star-shaped case in \cite{LLW1},
and all the theories of \cite{LoZ1} hold for dynamically convex star-shaped hypersurfaces by Theorem 1.4,
then the main results of \cite{WHL1} and \cite{Wan3} hold for dynamically convex star-shaped case,
i.e., Theorem 1.5 holds.\hfill\hb

\setcounter{equation}{0}
\section{Appendix }

In the section, we briefly review the equivariant Morse theory and the resonance identities for closed characteristics on
compact star-shaped hypersurfaces in $\R^{2n}$ developed in \cite{LLW1}. Now we fix a $\Sg\in\H_{st}(2n)$ and assume the following condition:

\noindent (F) {\bf There exist only finitely many geometrically distinct prime closed characteristics
\\$\quad \{(\tau_j, y_j)\}_{1\le j\le k}$ on $\Sigma$. }

Let $\hat{\sigma}=\inf_{1\leq j\leq k}{\sigma_j}$ and $T$ be a fixed positive constant.
Then by Section 2 of \cite{LLW1}, for any
$a>\frac{\hat{\sigma}}{T}$, we can construct a  function $\varphi_a\in C^{\infty}({\bf R}, {\bf R}^+)$ which has 0
as its unique critical point in $[0, +\infty)$. Moreover,
$\frac{\varphi^{\prime}(t)}{t}$ is strictly decreasing for $t>0$ together with $\varphi(0)=0=\varphi^{\prime}(0)$ and
$\varphi^{\prime\prime}(0)=1=\lim_{t\rightarrow
0^+}\frac{\varphi^{\prime}(t)}{t}$. More precisely, we
define $\varphi_a$ and the Hamiltonian function $\wtd{H}_a(x)=a\vf_a(j(x))$ via Lemma 2.2 and Lemma 2.4
in \cite{LLW1}. The precise dependence of $\varphi_a$ on $a$ is explained in Remark 2.3 of \cite{LLW1}.

For technical reasons we want to further modify the Hamiltonian, we define the new Hamiltonian
function $H_a$ via Proposition 2.5 of \cite{LLW1} and consider the fixed period problem
\be   \left\{\matrix{\dot{x}(t) &=& JH_a^\prime(x(t)), \cr
                   x(0) &=& x(T).\qquad         \cr }\right.  \lb{5.1}\ee
Then $H_a\in C^{3}({\bf R}^{2n} \setminus\{0\},{\bf R})\cap C^{1}({\bf R}^{2n},{\bf R})$.
Solutions of (\ref{5.1}) are $x\equiv 0$ and $x=\rho z(\sigma t/T)$ with
$\frac{\vf_a^\prime(\rho)}{\rho}=\frac{\sigma}{aT}$, where $(\sigma, z)$ is a solution of (\ref{1.1}). In particular,
non-zero solutions of (\ref{5.1}) are in one to one correspondence with solutions of (\ref{1.1}) with period
$\sigma<aT$.

For any $a>\frac{\hat{\sigma}}{T}$, we can choose some
large constant $K=K(a)$ such that
\be H_{a,K}(x) = H_a(x)+\frac{1}{2}K|x|^2   \lb{5.2}\ee
is a strictly convex function, that is,
\be (\nabla H_{a, K}(x)-\nabla H_{a, K}(y), x-y) \geq \frac{\ep}{2}|x-y|^2,  \lb{5.3}\ee
for all $x, y\in {\bf R}^{2n}$, and some positive $\ep$. Let $H_{a,K}^*$ be the Fenchel dual of $H_{a,K}$
defined by
\bea  H_{a,K}^\ast (y) = \sup\{x\cdot y-H_{a,K}(x)\;|\; x\in \R^{2n}\}.   \nn\eea
The dual action functional on $X=W^{1, 2}({\bf R}/{T {\bf Z}}, {\bf R}^{2n})$ is defined by
\be F_{a,K}(x) = \int_0^T{\left[\frac{1}{2}(J\dot{x}-K x,x)+H_{a,K}^*(-J\dot{x}+K x)\right]dt}.  \lb{5.4}\ee
Then $F_{a,K}\in C^{1,1}(X, \R)$ and for $KT\not\in 2\pi{\bf Z}$, $F_{a,K}$ satisfies the
Palais-Smale condition and $x$ is a
critical point of $F_{a, K}$ if and only if it is a solution of (\ref{5.1}). Moreover,
$F_{a, K}(x_a)<0$ and it is independent of $K$ for every critical point $x_a\neq 0$ of $F_{a, K}$.

When $KT\notin 2\pi{\bf Z}$, the map $x\mapsto -J\dot{x}+Kx$ is a Hilbert space isomorphism between
$X=W^{1, 2}({\bf R}/{T {\bf Z}}; {\bf R}^{2n})$ and $E=L^{2}({\bf R}/(T {\bf Z}),{\bf R}^{2n})$. We denote its inverse
by $M_K$ and the functional
\be \Psi_{a,K}(u)=\int_0^T{\left[-\frac{1}{2}(M_{K}u, u)+H_{a,K}^*(u)\right]dt}, \qquad \forall\,u\in E. \lb{5.5}\ee
Then $x\in X$ is a critical point of $F_{a,K}$ if and only if $u=-J\dot{x}+Kx$ is a critical point of $\Psi_{a, K}$.

Suppose $u$ is a nonzero critical point of $\Psi_{a, K}$.
Then the formal Hessian of $\Psi_{a, K}$ at $u$ is defined by
\be Q_{a,K}(v)=\int_0^T(-M_K v\cdot v+H_{a,K}^{*\prime\prime}(u)v\cdot v)dt,  \lb{5.6}\ee
which defines an orthogonal splitting $E=E_-\oplus E_0\oplus E_+$ of $E$ into negative, zero and positive subspaces.
The index and nullity of $u$ are defined by $i_K(u)=\dim E_-$ and $\nu_K(u)=\dim E_0$ respectively.
Similarly, we define the index and nullity of $x=M_Ku$ for $F_{a, K}$, we denote them by $i_K(x)$ and
$\nu_K(x)$. Then we have
\be  i_K(u)=i_K(x),\quad \nu_K(u)=\nu_K(x),  \lb{5.7}\ee
which follow from the definitions (\ref{5.4}) and (\ref{5.5}). The following important formula was proved in
Lemma 6.4 of \cite{Vit2}:
\be  i_K(x) = 2n([KT/{2\pi}]+1)+i^v(x) \equiv d(K)+i^v(x),   \lb{5.8}\ee
where the index $i^v(x)$ does not depend on K, but only on $H_a$.

By the proof of Proposition 2 of \cite{Vit1}, we have that $v\in E$ belongs to the null space of $Q_{a, K}$
if and only if $z=M_K v$ is a solution of the linearized system
\be  \dot{z}(t) = JH_a''(x(t))z(t).  \lb{5.9}\ee
Thus the nullity in (\ref{5.7}) is independent of $K$, which we denote by $\nu^v(x)\equiv \nu_K(u)= \nu_K(x)$.

By Proposition 2.11 of \cite{LLW1}, the index $i^v(x)$ and nullity $\nu^v(x)$ coincide with those defined for
the Hamiltonian $H(x)=j(x)^\alpha$ for all $x\in\R^{2n}$ and some $\aa\in (1,2)$. Especially
$1\le \nu^v(x)\le 2n-1$ always holds.

We have a natural $S^1$-action on $X$ or $E$ defined by
\be  \theta\cdot u(t)=u(\theta+t),\quad\forall\, \theta\in S^1, \, t\in\R.  \lb{5.10}\ee
Clearly both of $F_{a, K}$ and $\Psi_{a, K}$ are $S^1$-invariant. For any $\kappa\in\R$, we denote by
\bea
\Lambda_{a, K}^\kappa &=& \{u\in L^{2}({\bf R}/{T {\bf Z}}; {\bf R}^{2n})\;|\;\Psi_{a,K}(u)\le\kappa\},  \lb{5.11}\\
X_{a, K}^\kappa &=& \{x\in W^{1, 2}({\bf R}/(T {\bf Z}),{\bf R}^{2n})\;|\;F_{a, K}(x)\le\kappa\}.  \lb{5.12}\eea
For a critical point $u$ of $\Psi_{a, K}$ and the corresponding $x=M_K u$ of $F_{a, K}$, let
\bea
\Lm_{a,K}(u) &=& \Lm_{a,K}^{\Psi_{a, K}(u)}
   = \{w\in L^{2}(\R/(T\Z), \R^{2n}) \;|\; \Psi_{a, K}(w)\le\Psi_{a,K}(u)\},  \lb{5.13}\\
X_{a,K}(x) &=& X_{a,K}^{F_{a,K}(x)} = \{y\in W^{1, 2}(\R/(T\Z), \R^{2n}) \;|\; F_{a,K}(y)\le F_{a,K}(x)\}. \lb{5.14}\eea
Clearly, both sets are $S^1$-invariant. Denote by $\crit(\Psi_{a, K})$ the set of critical points of $\Psi_{a, K}$.
Because $\Psi_{a,K}$ is $S^1$-invariant, $S^1\cdot u$ becomes a critical orbit if $u\in \crit(\Psi_{a, K})$.
Note that by the condition (F), the number of critical orbits of $\Psi_{a, K}$
is finite. Hence as usual we can make the following definition.

{\bf Definition 5.1.} {\it Suppose $u$ is a nonzero critical point of $\Psi_{a, K}$, and $\Nn$ is an $S^1$-invariant open
neighborhood of $S^1\cdot u$ such that $\crit(\Psi_{a,K})\cap (\Lm_{a,K}(u)\cap \Nn) = S^1\cdot u$.
Then the $S^1$-critical
modules of $S^1\cdot u$ are defined by
\bea C_{S^1,\; q}(\Psi_{a, K}, \;S^1\cdot u)
=H_{q}((\Lambda_{a, K}(u)\cap\Nn)_{S^1},\; ((\Lambda_{a,K}(u)\setminus S^1\cdot u)\cap\Nn)_{S^1}). \nn\eea
Similarly, we define the $S^1$-critical modules $C_{S^1,\; q}(F_{a, K}, \;S^1\cdot x)$ of $S^1\cdot x$
for $F_{a, K}$.}

We fix $a$ and let $u_K\neq 0$ be a critical point of $\Psi_{a, K}$ with multiplicity $\mul(u_K)=m$,
that is, $u_K$ corresponds to a closed characteristic $(\tau, y)\subset\Sigma$ with $(\tau, y)$
being $m$-iteration of
some prime closed characteristic. Precisely, we have $u_K=-J\dot x+Kx$ with $x$
being a solution of (\ref{5.1}) and $x=\rho y(\frac{\tau t}{T})$ with
$\frac{\vf_a^\prime(\rho)}{\rho}=\frac{\tau}{aT}$.
Moreover, $(\tau, y)$ is a closed characteristic on $\Sigma$ with minimal period $\frac{\tau}{m}$.
By Lemma 2.10 of \cite{LLW1}, we construct a finite dimensional $S^1$-invariant subspace $G$ of
$L^{2}({\bf R}/{T {\bf Z}}; {\bf R}^{2n})$ and a functional $\psi_{a,K}$ on $G$.
For any $p\in\N$ satisfying $p\tau<aT$, we choose $K$
such that $pK\notin \frac{2\pi}{T}\Z$, then the $p$th iteration $u_{pK}^p$ of $u_K$ is given by $-J\dot x^p+pKx^p$,
where $x^p$ is the unique solution of (\ref{5.1}) corresponding to $(p\tau, y)$
and is a critical point of $F_{a, pK}$, that
is, $u_{pK}^p$ is the critical point of $\Psi_{a, pK}$ corresponding to $x^p$.
Denote by $g_{pK}^p$ the critical point of $\psi_{a,pK}$ corresponding to $u_{pK}^p$
and let $\wtd{\Lambda}_{a,K}(g_K)=\{g\in G \;|\; \psi_{a, K}(g)\le\psi_{a, K}(g_K)\}$.

Now we use the theory of Gromoll and Meyer, denote by $W(g_{pK}^p)$ the local characteristic manifold
of $g_{pK}^p$. Then we have

{\bf Proposition 5.2.}(cf. Proposition 4.2 of \cite{LLW1})
{\it For any $p\in\N$, we choose $K$ such that $pK\notin \frac{2\pi}{T}\Z$. Let $u_K\neq 0$
be a critical point of $\Psi_{a, K}$ with $\mul(u_K)=1$, $u_K=-J\dot x+Kx$ with $x$ being a critical point of
$F_{a, K}$. Then for all $q\in\Z$, we have
\bea
&& C_{S^1,\; q}(\Psi_{a,pK},\;S^1\cdot u_{pK}^p) \nn\\
&&\quad\cong \left(\frac{}{}H_{q-i_{pK}(u_{pK}^p)}(W(g_{pK}^p)\cap \wtd{\Lambda}_{a,pK}(g_{pK}^p),(W(g_{pK}^p)
                 \setminus\{g_{pK}^p\})\cap \wtd{\Lambda}_{a,pK}(g_{pK}^p))\right)^{\beta(x^p)\Z_p},
                  \qquad\quad \lb{5.15}\eea
where $\beta(x^p)=(-1)^{i_{pK}(u_{pK}^p)-i_K(u_K)}=(-1)^{i^v(x^p)-i^v(x)}$. Thus
\bea C_{S^1,\; q}(\Psi_{a,pK},\;S^1\cdot u_{pK}^p)=0 \quad
if q<i_{pK}(u_{pK}^p)~ or ~q>i_{pK}(u_{pK}^p)+\nu_{pK}(u_{pK}^p)-1.\lb{5.16}\eea
In particular, if $u_{pK}^p$ is
non-degenerate, i.e., $\nu_{pK}(u_{pK}^p)=1$, then}
\be C_{S^1,\; q}(\Psi_{a,pK},\;S^1\cdot u_{pK}^p)
    = \left\{\matrix{\Q, & {\rm if\;}q=i_{pK}(u_{pK}^p)\;{\rm and\;}\beta(x^p)=1,  \cr
                      0, & {\rm otherwise}. \cr}\right. \lb{5.17}\ee
We make the following definition:

{\bf Definition 5.3.} {\it For any $p\in\N$, we choose $K$ such that $pK\notin \frac{2\pi}{T}\Z$. Let $u_K\neq 0$ be
a critical point of $\Psi_{a,K}$ with $\mul(u_K)=1$, $u_K=-J\dot x+Kx$ with $x$ being a critical point of $F_{a, K}$.
Then for all $l\in\Z$, let
\bea
k_{l,\pm 1}(u_{pK}^p) &=& \dim\left(\frac{}{}H_{l}(W(g_{pK}^p)\cap
  \wtd{\Lambda}_{a,pK}(g_{pK}^p),(W(g_{pK}^p)\bs\{g_{pK}^p\})\cap \wtd{\Lambda}_{a,pK}(g_{pK}^p))\right)^{\pm\Z_p},
           \quad  \nn\\
k_l(u_{pK}^p) &=& \dim\left(\frac{}{}H_{l}(W(g_{pK}^p)\cap
  \wtd{\Lambda}_{a,pK}(g_{pK}^p),(W(g_{pK}^p)\bs\{g_{pK}^p\})\cap \wtd{\Lambda}_{a,pK}(g_{pK}^p))\right)^{\beta(x^p)\Z_p}.
           \qquad\quad   \nn\eea
Here $k_l(u_{pK}^p)$'s are called critical type numbers of $u_{pK}^p$. }

By Theorem 3.3 of \cite{LLW1}, we obtain that $k_l(u_{pK}^p)$ is independent of the choice of $K$ and
denote it by $k_l(x^p)$, here $k_l(x^p)$'s are called critical type numbers of $x^p$.

We have the following properties for critical type numbers:

{\bf Proposition 5.4.}(cf. Proposition 4.6 of \cite{LLW1})
{\it Let $x\neq 0$ be a critical point of $F_{a,K}$ with $\mul(x)=1$ corresponding to a
critical point $u_K$ of $\Psi_{a, K}$. Then there exists a minimal $K(x)\in \N$ such that
\bea
&& \nu^v(x^{p+K(x)})=\nu^v(x^p),\quad i^v(x^{p+K(x)})-i^v(x^p)\in 2\Z,  \qquad\forall p\in \N,  \lb{5.18}\\
&& k_l(x^{p+K(x)})=k_l(x^p), \qquad\forall p\in \N,\;l\in\Z. \lb{5.19}\eea
We call $K(x)$ the minimal period of critical modules of iterations of the functional $F_{a, K}$ at $x$. }

For every closed
characteristic $(\tau, y)$ on $\Sigma$, let $aT>\tau$ and choose $\vf_a$ as above.
Determine $\rho$ uniquely by $\frac{\vf_a'(\rho)}{\rho}=\frac{\tau}{aT}$. Let $x=\rho y(\frac{\tau t}{T})$.
Then we define the index $i(\tau, y)$ and nullity $\nu(\tau, y)$ of $(\tau, y)$ by
$$ i(\tau, y)=i^v(x), \qquad \nu(\tau, y)=\nu^v(x). $$
Then the mean index of $(\tau, y)$ is defined by
\be \hat i(\tau, y) = \lim_{m\rightarrow\infty}\frac{i(m\tau, y)}{m}.  \lb{5.20}\ee

Note that by Proposition 2.11 of \cite{LLW1}, the index and nullity are well defined and are independent of the
choice of $a$.

For a closed characteristic $(\tau, y)$ on $\Sigma$, we simply denote by $y^m\equiv(m\tau, y)$
the m-th iteration of $y$ for $m\in\N$.
By Proposition 3.2 of \cite{LLW1}, we can define the critical type numbers $k_l(y^m)$ of $y^m$ to be $k_l(x^m)$,
where $x^m$ is the critical point of $F_{a, K}$ corresponding to $y^m$. We also define $K(y)=K(x)$.

{\bf Remark 5.5.}(cf. Remark 4.10 of \cite{LLW1}) {\it
Note that $k_l(y^m)=0$ for $l\notin [0, \nu(y^m)-1]$ and it can take only values $0$
or $1$ when $l=0$ or $l=\nu(y^m)-1$. Moreover, the following facts are useful:

(i) $k_0(y^m)=1$ implies $k_l(y^m)=0$ for $1\le l\le \nu(y^m)-1$.

(ii) $k_{\nu(y^m)-1}(y^m)=1$ implies $k_l(y^m)=0$ for $0\le l\le \nu(y^m)-2$.

(iii) $k_l(y^m)\ge 1$ for some $1\le l\le \nu(y^m)-2$ implies $k_0(y^m)=k_{\nu(y^m)-1}(y^m)=0$.

(iv) In particular, only one of the $k_l(y^m)$s for $0\le l\le \nu(y^m)-1$ can be non-zero when $\nu(y^m)\le 3$.}

For a closed characteristic $(\tau,y)$ on $\Sigma$, the average Euler characteristic $\hat\chi(y)$ of $y$ is
defined by\bea \hat\chi(y)=\frac{1}{K(y)}\sum_{1\le m\le K(y)\atop 0\le l\le 2n-2}
(-1)^{i(y^{m})+l}k_l(y^{m}).  \lb{5.21}\eea
$\hat\chi(y)$ is a rational number. In particular, if all $y^m$s are
non-degenerate, then by Proposition 5.4 we have
\be \hat\chi(y)
    = \left\{\matrix{(-1)^{i(y)}, & {\rm if\;\;} i(y^2)-i(y)\in 2\Z,  \cr
           \frac{(-1)^{i(y)}}{2}, & {\rm otherwise}. \cr}\right.  \lb{5.22}\ee
We have the following mean index identities for closed characteristics.

{\bf Theorem 5.6.} {\it Suppose that $\Sg\in \H_{st}(2n)$ satisfies $\,^{\#}\T(\Sg)<+\infty$. Denote all
the geometrically distinct prime closed characteristics by $\{(\tau_j,\; y_j)\}_{1\le j\le k}$. Then the
following identities hold
\bea
\sum_{1\le j\le k\atop \hat{i}(y_j)>0}\frac{\hat{\chi}(y_j)}{\hat{i}(y_j)} &=& \frac{1}{2}, \lb{5.23} \\
\sum_{1\le j\le k\atop \hat{i}(y_j)<0}\frac{\hat{\chi}(y_j)}{\hat{i}(y_j)} &=& 0. \lb{5.24}\eea}

Let $F_{a, K}$ be a functional defined by (\ref{5.4}) for some $a, K\in\R$ sufficiently large and let $\epsilon>0$ be
small enough such that $[-\epsilon, 0)$ contains no critical values of $F_{a, K}$. For $b$ large enough,
The normalized Morse series of $F_{a, K}$ in $ X^{-\ep}\setminus X^{-b}$
is defined, as usual, by
\be  M_a(t)=\sum_{q\ge 0,\;1\le j\le p} \dim C_{S^1,\;q}(F_{a, K}, \;S^1\cdot v_j)t^{q-d(K)},  \lb{5.25}\ee
where we denote by $\{S^1\cdot v_1, \ldots, S^1\cdot v_p\}$ the critical orbits of $F_{a, K}$ with critical
values less than $-\epsilon$. The Poincar\'e series of
$H_{S^1, *}( X, X^{-\ep})$ is $t^{d(K)}Q_a(t)$, according to Theorem 5.1 of \cite{LLW1}, if we set
$Q_a(t)=\sum_{k\in \Z}{q_kt^k}$, then
\be   q_k=0 \qquad\qquad \forall\;k\in \mathring {I},  \lb{5.26}\ee
where $I$ is an interval of $\Z$ such that $I \cap [i(\tau, y), i(\tau, y)+\nu(\tau, y)-1]=\emptyset$ for all
closed characteristics $(\tau, y)$ on $\Sigma$ with $\tau\ge aT$. Then by Section 6 of \cite{LLW1}, we have
\be  M_a(t)-\frac{1}{1-t^2}+Q_a(t) = (1+t)U_a(t),   \lb{5.27}\ee
where $U_a(t)=\sum_{i\in \Z}{u_it^i}$ is a Laurent series with nonnegative coefficients.
If there is no closed characteristic with $\hat{i}=0$, then \be   M(t)-\frac{1}{1-t^2}=(1+t)U(t),   \lb{5.28}\ee
where $M(t)=\sum_{i\in\Z}{m_it^i}$ denotes the limit of $M_a(t)$ as $a$ tends to infinity,
$U(t)=\sum_{i\in\Z}{u_it^i}$ denotes the limit of $U_a(t)$ as $a$ tends to infinity and possesses
only non-negative coefficients. Specially, suppose that there exists an
integer $p<0$ such that the coefficients of $M(t)$ satisfy $m_p>0$ and $m_q=0$ for all integers $q<p$. Then
(\ref{5.28}) implies
\be   m_{p+1} \ge m_p.  \lb{5.29}\ee

\medskip

\bibliographystyle{abbrv}

\end{document}